\documentclass[a4paper,12pt,twoside]{article}
\usepackage{graphicx}\usepackage{a4wide}
\usepackage{amssymb,amsmath,amsthm}
\usepackage{siunitx}
\usepackage{subfigure,color,colordvi,graphicx,xcolor}
\usepackage[only,llbracket,rrbracket,interleave]{stmaryrd}
\usepackage[superscript,sort,compress]{cite}

\makeatletter
\newcommand*\wt[2][0.2ex]{%
        \begingroup
        \mathchoice{\wt@helper{#1}{#2}{\displaystyle}{\textfont}}
                   {\wt@helper{#1}{#2}{\textstyle}{\textfont}}
                   {\wt@helper{#1}{#2}{\scriptstyle}{\scriptfont}}
                   {\wt@helper{#1}{#2}{\scriptscriptstyle}{\scriptscriptfont}}%
        \endgroup
        #2%
}
\newcommand*\wt@helper[4]{%
        \def\currentfont{\the#41}%
        \def\currentskewchar{\char\the\skewchar\currentfont}%
        \setbox\tw@\hbox{\currentfont#2\currentskewchar}%
        \dimen@ii\wd\tw@
        \setbox\tw@\hbox{\currentfont#2{}\currentskewchar}%
        \advance\dimen@ii-\wd\tw@
        \rlap{\raisebox{-#1}{$\m@th#3\kern\dimen@ii\widetilde{\phantom{#2}}$}}%
}
\makeatother

\newcommand{\bm}[1]{\text{\boldmath $#1$\unboldmath}}
\newcommand{\norm}[1]{\lVert#1\rVert}

\newcommand{\mat}[1]{\mathbf{#1}}
\newcommand{\grad}{\bm{\nabla}}
\newcommand{\eltwo}{\ensuremath{\mathcal{L}_2}}

\newcommand{\ba}{\bm{a}}
\newcommand{\bu}{\bm{u}}

\newcommand{\bv}{\bm{v}}

\newcommand{\bn}{\bm{n}}
\newcommand{\bx}{\bm{x}}

\newcommand{\bmu}{\bm{\mu}}

\newcommand{\Insd}{\mat{I}_{d}}

\newcommand{\upgd}  {\bu_{_{\texttt{PGD}}}}
\newcommand{\ppgd}  {p_{_{\texttt{PGD}}}}
\newcommand{\uTpgd}  {\wt{\bu}_{_{\texttt{PGD}}}}
\newcommand{\pTpgd}  {\wt{p}_{_{\texttt{PGD}}}}
\newcommand{\uref}  {\bu_{_{\texttt{REF}}}}
\newcommand{\pref}  {p_{_{\texttt{REF}}}}
\newcommand{\fu}{\bm{f}_{\!\! u}}
\newcommand{\fd}{\bm{f}_{\!\! \diamond}}

\newcommand{\fp}{f_{\! p}}

\newcommand{\de}{\delta\!}
\newcommand{\I}{\mathcal{I}}
\newcommand{\bI}{\bm{\mathcal{I}}}
\newcommand{\De}{\varDelta}

\newenvironment{keywords}{\begin{quote}\emph{\textbf{Keywords:}}}{\end{quote}}

\theoremstyle{definition}

\newtheorem{remark}{Remark}

\RequirePackage{algorithm}[0.1] %Following SIAM cls
\usepackage{algorithmic} %Algorithm style, could alternatively use algpseudocode 

\graphicspath{{figsPDF/}{figsPNG/}}

\begin{document}
%==========================================================================
\title{Nonintrusive proper generalised decomposition for parametrised incompressible flow problems in OpenFOAM}

\author{
\renewcommand{\thefootnote}{\arabic{footnote}}
			  V. Tsiolakis\footnotemark[1]\textsuperscript{ \ ,}\footnotemark[2]\textsuperscript{ \ ,}\footnotemark[3] \ ,
			  M. Giacomini\footnotemark[2], \
			  R. Sevilla\footnotemark[3], \
             C. Othmer\footnotemark[1] \ and
             A. Huerta\footnotemark[2]
}

\date{\today}
%________________________________________________________________________
\maketitle

\renewcommand{\thefootnote}{\arabic{footnote}}

\footnotetext[1]{Volkswagen AG, Brieffach 011/1777, D-38436, Wolfsburg, Germany}
\footnotetext[2]{Laboratori de C\`alcul Num\`eric (LaC\`aN), ETS de Ingenieros de Caminos, Canales y Puertos, Universitat Polit\`ecnica de Catalunya, Barcelona, Spain}
\footnotetext[3]{Zienkiewicz Centre for Computational Engineering, College of Engineering, Swansea University, Wales, UK
\vspace{5pt}\\
Corresponding author: Matteo Giacomini. \textit{E-mail:} \texttt{matteo.giacomini@upc.edu}
}

%________________________________________________________________________
\begin{abstract}
The computational cost of parametric studies currently represents the major limitation to the application of simulation-based engineering techniques in a daily industrial environment. 
This work presents the first nonintrusive implementation of the proper generalised decomposition (PGD) in OpenFOAM, for the approximation of parametrised laminar incompressible Navier-Stokes equations.
The key feature of this approach is the seamless integration of a reduced order model (ROM) in the framework of an industrially validated computational fluid dynamics software.
This is of special importance in an industrial environment because in the online phase of the PGD ROM the description of the flow for a specific set of parameters is obtained simply via interpolation of the generalised solution, without the need of any extra solution step.
On the one hand, the spatial problems arising from the PGD separation of the unknowns are treated using the classical solution strategies of OpenFOAM, namely the semi-implicit method for pressure linked equations (SIMPLE) algorithm.
On the other hand, the parametric iteration is solved via a collocation approach.
The resulting ROM is applied to several benchmark tests of laminar incompressible Navier-Stokes flows, in two and three dimensions, with different parameters affecting the flow features.
Eventually, the capability of the proposed strategy to treat industrial problems is verified by applying the methodology to a parametrised flow control in a realistic geometry of interest for the automotive industry.
\end{abstract}

%________________________________________________________________________
\begin{keywords}
Reduced order models, proper generalised decomposition, finite volume, incompressible laminar Navier-Stokes, pressure Poisson equation, parametrised flows, OpenFOAM, nonintrusiveness
\end{keywords}

%==========================================================================
\section{Introduction and motivations}
\label{sc:Intro}
%==========================================================================

Computational fluid dynamics (CFD) is a key component in the current industrial design pipeline. 
Simulations of incompressible flows are performed on a daily basis to solve different problems both in automotive and aeronautical industries. 
Owing to its robustness, the most widely spread CFD methodology is the finite volume (FV) method~\cite{LeVeque:02,Toro:09,Sonar-MS:07,Ohlberger-BHO:17,Eymard-EGH:00,RS-SGH:18}.
Using this technique, numerically evaluated quantities of interest (e.g. drag and lift) have proved to match reasonably well experimental results. 

Nonetheless, design and optimisation cycles in a production environment require multiple queries of the same problem with boundary conditions, physical properties of the fluid and geometry of the domain varying within a range of values of interest.
In this context, parameters act as extra-coordinates of a high-dimensional partial differential equation (PDE).
The computational cost of such parametric studies currently represents the major limitation to the application of simulation-based engineering techniques in a daily industrial environment.
It is well-known that the computational complexity of approximating the PDEs describing the problems under analysis increases exponentially with the number of parameters considered. 
In recent years, reduced order models~\cite{AH-CHRW:17}, including reduced basis (RB)~\cite{Maday-BMNP:04,Grepl-GP:05,Grepl-GMNP:07,Rozza-RHP:08,DEIM-CS:10,Iapichino-IQR:12,Martini-MRH:15}, proper orthogonal decomposition (POD)~\cite{Iollo-ILD:00,Volkwein-KV:02,Iollo-BBI:09,John-CIJS:14,Hamdouni-AHLJ-14,Ballarin-BMQR:15,Hamdouni-LLPSH-18} and hierarchical model reduction (HiMod)~\cite{Perotto-PEV:10,Perotto-PV:14,Perotto-APV:18,Perotto-GPV:18}, have been proposed to reduce the computational burden of parametric analysis for several physical problems, including incompressible flows.
The aforementioned techniques rely on an \emph{a posteriori} reduction based on snapshots computed as solutions of the full-order model for different values of the parameters under analysis.
An alternative approach is represented by PGD~\cite{Ammar-AMCK:06,Ladeveze-CLC:10,Chinesta-CLBACGAAH:13,PGD-CCH:14,Chinesta-Keunings-Leygue}.
This method features an \emph{a priori} reduction~\cite{Ryckelynck-RCCA-06,Allery-VABHR-11,Allery-AHRV-11}, using a separable approximation of the solution, which depends explicitly on the parameters under analysis. 
In this context, during an offline phase, a reduced basis is constructed with no \emph{a priori} knowledge of the solution, whereas efficient online evaluations of the generalised solution are performed by simple interpolation in the parametric space.
The PGD framework has been first applied to incompressible Navier-Stokes equations in~\cite{Allery-DAA-10} to separate $x$ and $y$ directions in two-dimensional problems and in~\cite{Allery-DAA-11,Allery-LA-14} to separate space and time discretisations of unsteady flows.
See also~\cite{AH-AHCCL:14,SZ-ZDMH:15,Modesto-MZH:15,Signorini-SZD:17,PD-DZH:17,AH-HNC:18} for several applications of PGD to different physical problems.

In the context of flow problems, model reduction techniques based on Galerkin projection have been extensively studied in the literature~\cite{Karniadakis91,Karniadakis02,Zimmermann14}.
In this framework, several strategies have been proposed to construct the trial basis, using POD~\cite{HolmesLumleyBerkooz96}, RB~\cite{Haasdonk-HO:08} or the empirical interpolation method~\cite{Haasdonk-DHO:12}.
Concerning incompressible Navier-Stokes equations, in~\cite{Rozza-LCLR:16,Stabile-SR:18} supremiser stabilisations techniques have been investigated to couple the FV method with POD to solve parametrised turbulent flow problems.
Alternative projection methods based on minimisation of the residual of the momentum equation only~\cite{Allery-TALL-15} and on a least-squares Petrov-Galerkin approach~\cite{Carlberg11,Carlberg13,Carlberg17} have been proposed.
More recently, special attention has also been devoted to FV-based structure-preserving ROMs for conservation laws~\cite{Carlberg-CCS:18}.

Another key aspect for the application of simulation-based techniques to industrial problems is the capability of the proposed methods to provide verified and certified results.
This problem has been classically treated by equipping numerical methods with reliable and fully-computable \emph{a posteriori} error estimators using equilibrated fluxes~\cite{Destuynder-DM:99,Ern-ESV:10,Ern-EV:15} and flux-free approaches~\cite{NPM-PDH:06,Cottereau-CDH:09,NPM-PD:17} to control the error of the solution as well as of quantities of interest~\cite{Oden-OP:01,NPM-PBHP:06,NPM-PDH:09,PDM-LDH:10,Ainsworth-AR:12,Prudhomme-MP:15,MG-GPT:17,MG-G:18}.
Nonetheless, these approaches require intrusive modifications of existing computational libraries and may not be feasible in the context of commercial software. 
Hence, although the effort of the academic community in this direction, such solutions have not been successfully and widely integrated in codes utilised by the industry.
More recently, to circumvent this issue, great effort has been devoted to nonintrusive implementations in which novel numerical methodologies are externally coupled to existing commercial and open-source software used in industry on a daily basis.
Some contributions in this direction have been successfully proposed coupling PGD with Abaqus\textsuperscript{\textregistered} for mechanical problems \cite{Zou-ZCDA:18} and PGD with SAMTECH\textsuperscript{\textregistered} for shape optimisation problems \cite{Neron-CNLB:16}.
For flow problems, the coupling of POD and OpenFOAM has been discussed in \cite{Othmer-BOZ18,Stabile-SR:18}.
The present contribution is the first nonintrusive integration of the PGD framework in OpenFOAM for the solution of parametrised incompressible Navier-Stokes problems in the laminar regime.
The resulting algorithm, henceforth referred to as \texttt{pgdFoam}, relies on internal functions and routines of OpenFOAM~\cite{OpenFOAM} and exploits the incompressible flow solver \texttt{simpleFoam} for the spatial iteration of the alternating direction scheme.

The rest of this paper is organised as follows. 
Section~\ref{sc:FV-NS} recalls the incompressible Navier-Stokes equations and their FV approximation. 
The parametrised Navier-Stokes equations are introduced in Section~\ref{sc:Par-NS} as well as their PGD approximation and its nonintrusive implementation in OpenFOAM.
Numerical simulations to validate the discussed reduced-order strategy are provided in Section~\ref{sc:simulations}, whereas its application to parametrised flow control problems is presented in Section~\ref{sc:flowControl}.
Section~\ref{sc:Conclusion} summarises the results and two appendices complement the information with some technical details on the formulation and the OpenFOAM spatial solver utilised.

%==========================================================================
\section{Finite volume approximation of the incompressible Navier-Stokes equations}
\label{sc:FV-NS}
%==========================================================================

In this section, the steady Navier-Stokes equations for the simulation of incompressible viscous laminar flows in $d$ spatial dimensions are recalled. 
Let $\Omega\subset\mathbb{R}^{d}, \ \partial\Omega {=} \Gamma_D \cup \Gamma_N$ be an open bounded domain with disjoint Dirichlet, $\Gamma_D$, and Neumann, $\Gamma_N$, boundaries.
The flow problem under analysis consists of computing the velocity field $\bu$ and the pressure $p$ such that 
\begin{equation} \label{eq:NavierStokes}
\left\{\begin{aligned}
\grad {\cdot} (\bu {\otimes} \bu) - \grad {\cdot} (\nu \grad \bu) + \grad p &= \bm{s}       &&\text{in $\Omega$,}\\
\grad {\cdot} \bu &= 0  &&\text{in $\Omega$,}\\
\bu &= \bu_D  &&\text{on $\Gamma_D$,}\\
\bn {\cdot} (\nu \grad \bu {-} p \Insd ) &= \bm{t}  &&\text{on $\Gamma_N$,}
\end{aligned}\right.
\end{equation}
where the first equation describes the balance of momentum and the second one the conservation of mass.
In Equation~\eqref{eq:NavierStokes}, $\bm{s}$ represents a volumetric source term, $\nu {>} 0$ is the dynamic viscosity and $\Insd$ is the $d {\times} d$ identity matrix.
On the Dirichlet boundary $\Gamma_D$, the value $\bu_D$ of the velocity is imposed, whereas on $\Gamma_N$ the pseudo-traction $\bm{t}$ is applied.
From the modelling point of view, inlet surfaces and physical walls are described as Dirichlet boundaries with an imposed entering velocity profile and a homogeneous datum, respectively, whereas outlet surfaces feature homogeneous Neumann boundary conditions.
For the sake of simplicity and without loss of generality, $\Gamma_N$ is henceforth assumed to be an outlet boundary, that is a null $\bm{t}$ is considered.

%==========================================================================
\subsection{A cell-centred finite volume approximation using OpenFOAM}
\label{sc:CCFV}
%==========================================================================

In this section, the formulation of a FV scheme for the incompressible Navier-Stokes equations is briefly recalled to introduce the notation needed for the high-dimensional parametrised problem of Section~\ref{sc:Par-NS}.
The domain $\Omega$ is partitioned in $N$ nonoverlapping cells $V_i, \, i{=}1,\ldots,N$ such that $\Omega {:=} \bigcup_{i=1}^{N} V_i$ and $V_i {\cap} V_j {=} \emptyset, \, \text{for} \, i {\neq} j$.
%
%Moreover, the discrete interpolation space of piecewise constant element-by-element functions
%%
%\begin{equation*}
%\FV(\Omega) {:=} \{ v {\in} \eltwo(\Omega) : v\vert_{V_i} {=} \text{const} \ \forall V_i, \, i {=} 1, \dots, N \}
%\end{equation*} 
%%
%is introduced.
%%
%The finite volume discretisation is constructed starting from the weak formulation of Equation~\eqref{eq:NavierStokes}, namely find $(\bu,p) {\in} [\FV(\Omega)]^{d} {\times} \FV(\Omega)$ such that $\bu {=} \bu_D \, \text{on} \, \Gamma_D$ and on each cell $V_i, \,  i=1,\dots,N$ it holds
%%
%\begin{equation} \label{eq:weak-NS}
%\!\left\{\begin{aligned}
%\int_{V_i}{\!\!\! \bv {\cdot} [\grad {\cdot} (\bu {\otimes} \bu)] \, dV} 
%- \int_{V_i}{\!\!\! \bv {\cdot} [\grad {\cdot} (\nu \grad \bu)]  \, dV} 
%+ \int_{V_i}{\!\!\! \bv {\cdot} \grad p \, dV} 
%&= \int_{V_i}{\!\!\! \bv {\cdot} \bm{s} \, dV} , \\
%\int_{V_i}{\!\!\! q \, \grad {\cdot} \bu \, dV} &= 0 ,
%\end{aligned}\right.
%\end{equation}
%%
%for all $(\bv,q) {\in} [\FV(\Omega)]^{d} {\times} \FV(\Omega)$.
%%
%Note that considering the unitary test functions in Equation~\eqref{eq:weak-NS}, the integral formulation of the Navier-Stokes equations is retrieved according to the classical finite volume rationale \cite{Ohlberger-BHO:17,Eymard-EGH:00}.
%
The FV discretisation is constructed starting from the integral formulation of Equation~\eqref{eq:NavierStokes}, namely find $(\bu,p)$, constant on each cell $V_i, \,  i{=}1,\dots,N$, such that $\bu {=} \bu_D \, \text{on} \, \Gamma_D$ and it holds
\begin{equation} \label{eq:weak-NS}
\!\left\{\begin{aligned}
\int_{V_i}{\!\! \grad {\cdot} (\bu {\otimes} \bu) \, dV} 
- \int_{V_i}{\!\!  \grad {\cdot} (\nu \grad \bu) \, dV} 
+ \int_{V_i}{\!\! \grad p \, dV} 
&= \int_{V_i}{\!\! \bm{s} \, dV} , \\
\int_{V_i}{\!\! \grad {\cdot} \bu \, dV} &= 0 .
\end{aligned}\right.
\end{equation}

OpenFOAM implements a cell-centred FV rationale in which piecewise constant approximations are sought for velocity and pressure in each cell of the computational mesh and the degrees of freedom of the discretised problem are located at the centroid of each finite volume. 
Employing Gauss's theorem, the integrals in Equation~\eqref{eq:weak-NS} are rewritten in terms of fluxes over the boundaries of the cells and approximated using classical central differencing schemes~\cite{Ohlberger-BHO:17,Eymard-EGH:00}.
Moreover, to handle the nonlinearity in the convection term, OpenFOAM considers a relaxation approach introducing a fictitious time variable.
The resulting solution strategy relies on the SIMPLE algorithm which belongs to the family of fractional-step projection methods~\cite{Patankar-PS:72,Donea-Huerta}. A brief description of this method is provided in~\ref{sc:SIMPLE}.

%==========================================================================
\section{Nonintrusive proper generalised decomposition for parametrised flow problems}
\label{sc:Par-NS}
%==========================================================================

Consider now the case in which the user-prescribed data in Equation~\eqref{eq:NavierStokes}, i.e. the viscosity coefficient, the source term and the boundary conditions, depend on a set of parameters $\bmu\in\bI\subset\mathbb{R}^{M}$, with $M$ being the number of parameters. More presicely, the set $\bI$ describing the range of admissible parameters can be defined as the Cartesian product of the domains of the $M$ parameters, namely,  $\bI {:=} \I_1\times\I_2\times\dotsb\times\I_{M}$ with $\mu_i\in\I_i$ for $i {=} 1,\dots , M$. 
Within this context, $\bmu$ is treated as a set of additional independent variables, or parametric coordinates, instead of problem parameters. 
For the purpose of discretisation, each interval $\I_i$ is subdivided in $N_\mu$ subintervals. 
%and the discrete parametric interpolation space 
%%
%\begin{equation*}
%\Lh(\bI) {:=} \{ (\phi_1, \dots, \phi_{M})^T : \phi_i {\in} \eltwo(\I_i) \, \text{and} \, \phi_i\vert_{\I_{i,j}} {=} \text{const} \ \forall \I_{i,j}, \, j {=} 1, \dots, N_\mu \}
%\end{equation*} 
%%
%is defined.
%%
%The unknown pair $(\bu,p) \in [\FV(\Omega)\otimes\Lh(\bI)]^{d} {\times} [\FV(\Omega)\otimes\Lh(\bI)]$ is thus sought in a high-dimensional space described by the independent variables $(\bx , \bmu ) \in\Omega {\times} \bI$ and fulfils the following parametrised Navier-Stokes equations
%%
%\begin{equation} \label{eq:weak-NS-PGD}
%\!\left\{\begin{aligned}
%\int_{\bI}\int_{V_i}{\!\!\! \bv {\cdot} \bigl[\grad {\cdot} (\bu {\otimes} \bu)\bigr] \, dV \, d\bI} 
%- \int_{\bI}\int_{V_i} \!\!\! \bv {\cdot} \bigl[\grad {\cdot} (\nu \grad \bu)\bigr] \, & dV \, d\bI \\
%+ \int_{\bI}\int_{V_i}{\!\!\! \bv {\cdot} \grad p  \, dV \, d\bI} 
%&= \int_{\bI}\int_{V_i}{\!\!\! \bv {\cdot}\bm{s} \, dV \, d\bI} , \\
%\int_{\bI}\int_{V_i}{\!\!\! q \, \grad {\cdot} \bu \, dV \, d\bI} &= 0 ,
%\end{aligned}\right.
%\end{equation}
%%
%for all $(\bv,q) {\in} [\FV(\Omega)\otimes\Lh(\bI)]^{d} {\times} [\FV(\Omega)\otimes\Lh(\bI)]$.
%%
%
The unknown pair $(\bu,p)$ is thus sought in a high-dimensional space described by the independent variables $(\bx , \bmu ) \in\Omega {\times} \bI$ and fulfils the following parametrised Navier-Stokes equations on each cell $V_i, \, i{=}1,\ldots,N$
\begin{equation} \label{eq:weak-NS-PGD}
\!\left\{\begin{aligned}
\int_{\bI}\int_{V_i}{\!\! \grad {\cdot} (\bu {\otimes} \bu) \, dV \, d\bI} 
- \int_{\bI}\int_{V_i} \!\! \grad {\cdot} (\nu \grad \bu) \, dV \, d\bI &\\
+ \int_{\bI}\int_{V_i}{\!\! \grad p  \, dV \, d\bI} 
&= \int_{\bI}\int_{V_i}{\!\! \bm{s} \, dV \, d\bI} , \\
\int_{\bI}\int_{V_i}{\!\! \grad {\cdot} \bu \, dV \, d\bI} &= 0 .
\end{aligned}\right.
\end{equation}
In the following sections, the rationale for the construction of a separated solution of the parametrised Navier-Stokes equations is recalled and the proposed nonintrusive implementation of the alternating direction scheme in OpenFOAM is presented.

%==========================================================================
\subsection{The proper generalised decomposition rationale}
\label{sc:PGD}
%==========================================================================

PGD constructs an approximation $(\upgd^n, \ppgd^n)$ of the solution $(\bu, p)$ of Equation~\eqref{eq:weak-NS-PGD} in terms of a sum of $n$ separable functions, or modes.
Each mode is the product of functions depending solely on one of the arguments $\bx, \mu_1,\dots,\mu_{M}$. For the sake of readability and without loss of generality, only the spatial coordinates $\bx$ and the parametric ones $\bmu$ are henceforth separated.

Following \cite{PD-DZH:17}, the so-called \emph{single parameter approximation} is detailed. That is, for each mode, a unique scalar parametric function $\phi(\bmu)$ is considered for all the variables and the resulting separated form of the unknowns is
\begin{equation}\label{eq:sep}
 \left\{\begin{aligned}
 \upgd^n(\bx,\bmu)   &= \upgd^{n-1}(\bx,\bmu)   + \sigma_u^n  \fu^n(\bx) \phi^n(\bmu), \\
 \ppgd^n(\bx,\bmu)   &= \ppgd^{n-1}(\bx,\bmu)   + \sigma_p^n  \fp^n(\bx) \phi^n(\bmu),
 \end{aligned}\right.
\end{equation}
where the superindex $n$ denotes the, \emph{a priori} unknown, number of terms in the PGD expansion and the positive scalar coefficients $\sigma_u^n$ and $\sigma_p^n$ represent the amplitude of the $n$-th mode for velocity and pressure, respectively. 
These coefficients are obtained normalising the modal functions, namely
\begin{equation*}
\sigma_u^n := \| \fu^n \| \,\, \text{and} \,\, \sigma_p^n := \| \fp^n \| , 
\end{equation*}
with $\| \phi^n \| = 1$. 
Appropriate user-defined norms on the spatial and parametric domains need to be introduced for each function. For all the simulations in Section~\ref{sc:simulations} and \ref{sc:flowControl}, the $\eltwo$ norm has been considered for normalisation.
\begin{remark}
The normalisation coefficients play a critical role in checking the convergence of the PGD algorithm and may be used as quantitative stopping criterion in the PGD enrichment procedure described in Section~\ref{sc:PGD-ADI}.
\end{remark}
For a discussion on alternative formulations of the separation in Equation~\eqref{eq:sep}, involving both scalar and vector-valued parametric functions, the interested reader is referred to~\cite{PD-DZH:17}.
Henceforth and except in case of ambiguity, the dependence of the modes on $\bx$ and $\bmu$ is omitted.

Considering a linearised approach to compute each new mode, Equation~\eqref{eq:sep} can be rewritten as the following \emph{predictor-corrector single parameter approximation}
\begin{equation}\label{eq:sep-increment}
\left\{\begin{alignedat}{2}
 \upgd^n   &= \uTpgd^n + \sigma_u^n\de\uTpgd^n              &&= \upgd^{n-1} + \sigma_u^n\fu^n\phi^n + \sigma_u^n\de\uTpgd^n , \\
 \ppgd^n   &= \pTpgd^{\:n} + \sigma_p^n\de\pTpgd^{\:n}   &&= \ppgd^{n-1} +\sigma_p^n\fp^n\phi^n + \sigma_p^n\de\pTpgd^{\:n} ,
 \end{alignedat}\right.
\end{equation}
where $\uTpgd^n {:=} \upgd^{n-1} + \sigma_u^n\fu^n\phi^n$ and $\pTpgd^{\:n} {:=} \ppgd^{n-1} +\sigma_p^n\fp^n\phi^n$ account for the $n{-}1$ previously computed terms and a prediction of the current mode.
More precisely, $(\sigma_u^n\fu^n\phi^n,\sigma_p^n\fp^n\phi^n)$ play the role of predictors in the computation of the $n$-th mode, whereas $(\sigma_u^n\de\uTpgd^n,\sigma_p^n\de\pTpgd^{\:n})$ are the corresponding correctors featuring the variations $\De$ in the spatial and parametric functions, namely
\begin{equation}\label{eq:PGD-corr}
\left\{\begin{aligned}
\de\uTpgd^n &:= \De\fu\phi^n+\fu^n\De\phi+\De\fu\De\phi , \\
\de\pTpgd^{\:n} &:= \De\fp\phi^n+\fp^n\De\phi+\De\fp\De\phi .
\end{aligned}\right.
\end{equation}
Note that the last term in Equation~\eqref{eq:PGD-corr} represents a high-order variation which is henceforth neglected.
As for the classical \emph{single parameter approximation}, $\sigma_u^n$ and $\sigma_p^n$ represent the amplitudes of the $n$-th velocity and pressure modes. That is, setting $\| \phi^n + \De\phi \| = 1$, they are defined as
\begin{equation*}
\sigma_u^n := \| \fu^n + \De\fu \| \,\, \text{and} \,\, \sigma_p^n := \| \fp^n + \De\fp \| .
\end{equation*}

%==========================================================================
\subsection{Predictor-corrector alternating direction scheme}
\label{sc:PGD-ADI}
%==========================================================================

In order to compute $(\upgd^n, \ppgd^n)$ in Equation~\eqref{eq:sep-increment}, a greedy algorithm is implemented. The first PGD mode $(\upgd^0, \ppgd^0)$ is arbitrarily chosen to fulfil the Dirichlet boundary conditions of the problem and the $n$-th mode is successively computed assuming that term $n{-}1$ is available \cite{PGD-CCH:14, Chinesta-Keunings-Leygue}. Some variations of this strategy based on Arnoldi-type iterations have been investigated in \cite{Nouy:08,LT-TMN:14}. 
In this section, the alternating direction scheme used to compute the PGD modes is described. A key assumption for the application of this method is the separability of the data. For the sake of simplicity and without any loss of generality, the separated form of the viscosity coefficient, see e.g. \cite{SZ-ZDMH:15}, is reported
\begin{equation}\label{eq:Ksep}
 \nu(\bx,\bmu) := \psi(\bmu)D(\bx) = \sum_{i=1}^{n_{\nu}} \psi_{1,i}(\mu_1)\dotsb\psi_{M,i}(\mu_{M})D_i(\bx) ,
\end{equation}
and analogous separations are considered for all the parametric data in the problem under analysis.

By plugging \eqref{eq:sep-increment} into \eqref{eq:weak-NS-PGD} and gathering the unknown increments  $(\sigma_u^n\de\uTpgd^n,\sigma_p^n\de\pTpgd^{\:n})$ on the left-hand side while leaving on the right-hand side the residuals computed using the previous modes $(\upgd^{n-1},\ppgd^{n-1})$ and the predictions $(\sigma_u^n\fu^n\phi^n,\sigma_p^n\fp^n\phi^n)$ of the current one, the following equations are obtained
\begin{equation} \label{eq:PGD-increment}
\left\{\begin{aligned}
\int_{\bI}\int_{V_i}{\!\! \grad {\cdot} (\sigma_u^n\de\uTpgd^n {\otimes} \sigma_u^n\de\uTpgd^n) \, dV \, d\bI} 
\hspace{182pt} & \\
+ \int_{\bI}\int_{V_i}{\!\! \grad {\cdot} (\sigma_u^n\de\uTpgd^n {\otimes} \uTpgd^n) \, dV \, d\bI} 
+ \int_{\bI}\int_{V_i}{\!\! \grad {\cdot} (\uTpgd^n {\otimes} \sigma_u^n\de\uTpgd^n) \, dV \, d\bI} 
\hspace{18pt} & \\
- \int_{\bI}{ \psi \int_{V_i} \!\! \grad {\cdot} (D \grad (\sigma_u^n\de\uTpgd^n))  \, dV \, d\bI} 
+ \int_{\bI}\int_{V_i}{\!\! \grad (\sigma_p^n\de\pTpgd^{\:n})  \, dV \, d\bI} = \mathcal{R}_u , 
 & \\
\int_{\bI}\int_{V_i}{\!\! \grad {\cdot} (\sigma_u^n\de\uTpgd^n) \, dV \, d\bI} = \mathcal{R}_p , 
&
%\hspace{130pt} &
\end{aligned}\right.
\end{equation}
where the residuals are defined as 
\begin{align}
&
\begin{aligned}
\mathcal{R}_u :=& \, \mathcal{R}_u(\upgd^{n-1}, \ppgd^{n-1}, \sigma_u^n\fu^n, \sigma_p^n\fp^n, \phi^n) = \mathcal{R}_u(\uTpgd^n,\pTpgd^{\:n})\\
=& \int_{\bI}\int_{V_i}{\!\! \bm{s} \, dV \, d\bI}
- \int_{\bI}\int_{V_i}{\!\! \grad {\cdot} (\uTpgd^n {\otimes} \uTpgd^n) \, dV \, d\bI} \\
&+ \int_{\bI} \psi \int_{V_i}{\!\! \grad {\cdot} (D \grad \uTpgd^n)  \, dV \, d\bI} 
- \int_{\bI}\int_{V_i}{\!\! \grad \pTpgd^{\:n}  \, dV \, d\bI} ,
\end{aligned}
\label{eq:resU} 
\\ &
\begin{aligned}
\mathcal{R}_p :=& \, \mathcal{R}_p(\upgd^{n-1}, \sigma_u^n\fu^n, \phi^n) = \mathcal{R}_p(\uTpgd^n) \\
=& - \int_{\bI}\int_{V_i}{\!\! \grad {\cdot} \uTpgd^n \, dV \, d\bI} .
\end{aligned}
\label{eq:resP}
\end{align}
%

%The test functions in Equation~\eqref{eq:PGD-increment} characterise the tangent manifold associated with the spatial, $\bx$, and parametric, $\bmu$, coordinates. 
%%
%Following from Equation~\eqref{eq:PGD-corr}, the PGD aims to compute the pair $(\sigma_u^n \De\fu, \sigma_p^n \De\fp) \in [\FV(\Omega)]^{d} {\times} \FV(\Omega)$ and $\De\phi \in \Lh(\bI)$ such that Equation~\eqref{eq:PGD-increment} is verified on the tangent manifold of the predicted mode 
%%
%\begin{equation*}
%\bv = \de\fu\phi^n {+} \sigma_u^n \fu^n\de\phi
%\,\, \text{ and } \,\, 
%q = \de\fp\phi^n {+} \sigma_p^n \fp^n\de\phi
%\end{equation*}
%for all $(\de\fu,\de\fp) \in [\FV(\Omega)]^{d} {\times} \FV(\Omega)$ and $\de\phi \in \Lh(\bI)$.
%%
%Exploiting the separated structure of the unknown and test functions as well as the affine parametric decomposition of the involved integral forms, the numerical complexity of the high-dimensional equation is reduced by alternatively solving for the spatial and the parametric unknowns, as detailed in the next subsections.
%

As classical in ROMs \cite{Patera-Rozza:07,Rozza:14}, an affine dependence of the forms in \eqref{eq:PGD-increment}, \eqref{eq:resU} and \eqref{eq:resP} on the parameters is required to construct the PGD approximation.
The spatial (respectively, parametric) component of each mode is thus computed by restricting Equation~\eqref{eq:PGD-increment} to the tangent manifold associated with the spatial (respectively, parametric) coordinate. 
Following from Equation~\eqref{eq:PGD-corr} and setting a fixed value for the parametric function $\phi^n$, the pair $(\sigma_u^n \De\fu, \sigma_p^n \De\fp)$ is determined by solving a purely spatial PDE.
Recall that the PGD alternating direction scheme handles homogeneous Dirichlet boundary conditions at each iteration of the spatial solver \cite{Chinesta-Keunings-Leygue}, whereas inhomogeneous data are treated by the first arbitrary PGD mode $(\upgd^0, \ppgd^0)$ introduced above.
In a similar fashion, the increment $\De\phi$ is computed as the solution of an algebraic system of equations in the parameter $\bmu$ while the spatial functions $(\sigma_u^n \fu^n,\sigma_p^n \fp^n)$ are considered known.

Note that at each iteration of the alternating direction scheme, $\uTpgd^n$ is known and may be expressed in separated form as $\sum_{m=1}^n \sigma_u^m\fu^m\phi^m$.
Thus, exploiting the separated structure of the unknowns and the affine parametric decomposition of the involved integral forms, the numerical complexity of the high-dimensional PDE is reduced by alternatively solving for the spatial and the parametric unknowns, as detailed in the next subsections.
\begin{remark}
By restricting Equation~\eqref{eq:PGD-increment} to the tangent manifold in the spatial (respectively, parametric) direction, the integral forms are multiplied by $\phi^n$ (respectively, $(\sigma_u^n\fu^n,\sigma_p^n\fp^n)$).
This is equivalent to the projection of the high-dimensional PDE to the tangent manifold discussed for PGD in the context of finite element approximations~\cite{PD-DZH:17}.
More precisely, being $\bv(\bx,\bmu)$ (respectively, $q(\bx,\bmu)$) the test function in the finite element weak form of the momentum (respectively, continuity) equation, the projection on the tangent manifold leads to
\begin{equation*}
%\begin{aligned}
%\bv(\bx,\bmu) &= \de\fu(\bx)\phi^n(\bmu) {+} \sigma_u^n \fu^n(\bx)\de\phi(\bmu) ,
%\\
%q(\bx,\bmu) &= \de\fp(\bx)\phi^n(\bmu) {+} \sigma_p^n \fp^n(\bx)\de\phi(\bmu) ,
%\end{aligned}
(\bv,q) =
\begin{cases}
(\phi^n\de\fu,\phi^n\de\fp) , &  \text{for the spatial iteration} ,\\
(\sigma_u^n \fu^n\de\phi,\sigma_p^n \fp^n\de\phi) , & \text{for the parametric iteration} ,
\end{cases}
\end{equation*}
where $(\de\fu,\de\fp)$ and $\de\phi$ are test functions depending solely on the spatial and parametric variables, respectively.
In the framework of FV discretisations, these test functions are set equal to $1$ to retrieve the classical integral formulation of the PDE under analysis.
The corresponding restriction to the tangent manifold thus leads to the following functions multiplying the integral forms in Equations~\eqref{eq:PGD-increment}, \eqref{eq:resU} and \eqref{eq:resP}
\begin{equation*}
(\bv,q) =
\begin{cases}
(\phi^n,\phi^n) , & \text{for the spatial iteration} ,\\
(\sigma_u^n \fu^n,\sigma_p^n \fp^n) , & \text{for the parametric iteration} .
\end{cases}
\end{equation*}
\end{remark}

%==========================================================================
\subsubsection{The spatial iteration}
\label{sc:PGD-spatial}
%==========================================================================

First, the parametric function $\phi^n$ is fixed and the increments $(\sigma_u^n \De\fu, \sigma_p^n \De\fp)$ are determined by solving a spatial PDE. 
More precisely, restricting Equation~\eqref{eq:PGD-increment} to the tangent manifold in the spatial direction and neglecting the high-order terms, a pair $(\sigma_u^n \De\fu, \sigma_p^n \De\fp)$, constant element-by-element, is sought such that in each cell $V_i, \, i{=}1,\ldots,N$ it holds
\begin{equation} \label{eq:PGD-spatial}
\hspace{-5pt}
\left\{\begin{aligned}
% \sum_{m=1}^n \! \alpha_1^m \! \left( 
%		\int_{V_i}{\!\! \grad {\cdot} ( \sigma_u^m \fu^m {\otimes} \sigma_u^n \De\fu) \, dV} 
%		 {+} \!\! \int_{V_i}{\!\! \grad {\cdot} (\sigma_u^n \De\fu {\otimes} \sigma_u^m \fu^m) \, dV} \right) 
\hspace{-3pt}
		\int_{V_i}{\!\! \grad {\cdot} \Big( \sigma_u^n \De\fu {\otimes} \sum_{m=1}^n \! \alpha_1^m \sigma_u^m \fu^m \Big) \, dV}
		 {+} \!\! \int_{V_i}{\!\! \grad {\cdot} \Big( \sum_{m=1}^n \! \alpha_1^m \sigma_u^m \fu^m {\otimes} \sigma_u^n \De\fu \Big) \, dV} 
& \\
- \alpha_2 \!\! \int_{V_i} \!\! \grad {\cdot} (D \grad (\sigma_u^n \De\fu))  \, dV
 + \! \alpha_3 \!\! \int_{V_i}{\!\! \grad (\sigma_p^n \De\fp)  \, dV} = \, R_u^n , &
 \\
\alpha_3 \!\! \int_{V_i}{\grad {\cdot} (\sigma_u^n \De\fu) \, dV} = \, R_p^n , &
\end{aligned}\right.
\end{equation}
where $R_u^n$ and $R_p^n$ are the spatial residuals associated with the discretisation of the momentum and mass equations, respectively, and each coefficient $\alpha_i, \, i{=}1,\ldots,3$ depends solely on the parametric function $\phi^n$ and on the data of the problem
\begin{equation}\label{eq:coef-spatial}
  \alpha_1^m := \int_{\bI}{ \left[\phi^n\right]^2 \phi^m \, d\bI} , \quad
  \alpha_2 := \int_{\bI}{ \left[\phi^n\right]^2 \psi \, d\bI} , \quad
  \alpha_3 := \int_{\bI}{ \left[\phi^n\right]^2 \, d\bI} .
\end{equation}
Note that given the separable form of \eqref{eq:resU}-\eqref{eq:resP}, an efficient implementation of the right-hand side of the spatial iteration may be devised and the corresponding FV discretisation is obtained.
A detailed description of the residuals acting as linear functionals on the right-hand side of Equation~\eqref{eq:PGD-spatial} is provided in \ref{sc:residuals}.

%The discretisation of Equation~\eqref{eq:PGD-spatial} is thus performed using the cell-centred finite volume method implemented in OpenFOAM, see Section~\ref{sc:CCFV}.
%%
%Note that the first term in Equation~\eqref{eq:PGD-spatial} represents a high-order variation of the velocity field and introduces a negligible correction of the predictor $\fu^n$, which is thus omitted.
%
%The remaining terms feature the same structure as the original incompressible Navier-Stokes problem in the spatial domain $\Omega$, see Equation~\eqref{eq:weak-NS}.
%
The terms in Equation~\eqref{eq:PGD-spatial} feature a structure similar to the original incompressible Navier-Stokes problem in the spatial domain $\Omega$, see Equation~\eqref{eq:weak-NS}.
The discretisation is thus performed using the cell-centred FV method implemented in OpenFOAM, see Section~\ref{sc:CCFV}.
The main difference is represented by the first two integrals on the left-hand side of the momentum equation in~\eqref{eq:PGD-spatial}. 

On the one hand, the first integral is a relaxation of the classical nonlinear convection term in Navier-Stokes equations, based on the last computed approximation $\sum_{m=1}^n \! \alpha_1^m \sigma_u^m \fu^m$ of the unknown velocity field.
From a practical point of view, this treatment of the convection term is equivalent to the one performed by the SIMPLE algorithm which solves a linearised version of the Navier-Stokes equations by introducing a fictitious time variable $t^k$ and substituting the unknown convection field with its approximation at time $t^{k-1}$, see~\ref{sc:SIMPLE}.

On the other hand, the second integral does not appear in the classical Navier-Stokes equations and is thus not treated by the SIMPLE algorithm.
In order to preserve the nonintrusiveness of the discussed PGD approach, the standard solution strategy implemented in OpenFOAM for such problem, namely \texttt{simpleFoam}, is applied.
Hence, a relaxation is introduced in the SIMPLE iterations and this term is handled explicitly as part of the right-hand side of the momentum equation, leading to
%
%\begin{equation} \label{eq:PGD-spatial-final}
%\left\{\begin{aligned}
%\sum_{m=1}^n \alpha_2^m \!\! \int_{V_i}{\!\!\! \de\fu {\cdot} \bigl[\grad {\cdot} (\sigma_u^n \De\fu {\otimes} \sigma_u^m \fu^m)\bigr] \, dV}
%\hspace{130pt} & \\
%- \alpha_3 \!\! \int_{V_i} \!\!\!  \de\fu {\cdot} [\grad {\cdot} (D \grad (\sigma_u^n \De\fu))]  \, dV
% + \alpha_4 \!\! \int_{V_i}{\!\!\! \de\fu {\cdot} \grad (\sigma_p^n \De\fp)  \, dV} 
%\hspace{20pt} & \\
% = \mathcal{R}_u(\phi^n \de\fu) 
% - \sum_{m=1}^n \alpha_2^m \!\! \int_{V_i}{\!\!\! \de\fu {\cdot} \bigl[\grad {\cdot} ( \sigma_u^m \fu^m {\otimes} \sigma_u^k \De\fu^k)\bigr] \, dV} , 
%& \\
%\alpha_4 \!\! \int_{V_i}{\!\!\! \de\fp \, \grad {\cdot} (\sigma_u^n \De\fu) \, dV} = \mathcal{R}_p(\phi^n \de\fp) ,
% \hspace{140pt} &
%\end{aligned}\right.
%\end{equation}
%
\begin{equation} \label{eq:PGD-spatial-final}
\left\{\begin{aligned}
%\sum_{m=1}^n \alpha_1^m \!\! \int_{V_i}{\grad {\cdot} (\sigma_u^n \De\fu {\otimes} \sigma_u^m \fu^m) \, dV}
\int_{V_i}{\!\! \grad {\cdot} \Big( \sigma_u^n \De\fu {\otimes}\sum_{m=1}^n \! \alpha_1^m \sigma_u^m \fu^m \Big) \, dV} 
- \alpha_2 \!\! \int_{V_i} \grad {\cdot} (D \grad (\sigma_u^n \De\fu)) \, dV 
\hspace{25pt} & \\
 + \alpha_3 \!\! \int_{V_i}{\grad (\sigma_p^n \De\fp)  \, dV} 
  = R_u^n 
% - \sum_{m=1}^n \alpha_1^m \!\! \int_{V_i}{\grad {\cdot} ( \sigma_u^m \fu^m {\otimes} \sigma_u^k \De\fu^k) \, dV} , 
{-} \!\! \int_{V_i}{\!\! \grad {\cdot} \Big( \sum_{m=1}^n \! \alpha_1^m \sigma_u^m \fu^m {\otimes} \sigma_u^{k-1} \De\fu^{k-1} \Big) \, dV} 
& \\
\alpha_3 \!\! \int_{V_i}{\grad {\cdot} (\sigma_u^n \De\fu) \, dV} = R_p^n ,
\hspace{198pt} &
\end{aligned}\right.
\end{equation}
where the index $k{-}1$ is now associated with the lastly computed increment $\sigma_u^{k-1}\De\fu^{k-1}$ in the SIMPLE algorithm, see~\ref{sc:SIMPLE}.

%==========================================================================
\subsubsection{The parametric iteration}
\label{sc:PGD-param}
%==========================================================================
 
In the parametric step, the value of the previously computed spatial functions is fixed $(\fu^n,\fp^n) {\gets} (\sigma_u^n \fu^n {+} \De\fu, \sigma_p^n \fp^n {+} \De\fp)$ and the parametric increment $\De\phi$ acts as unknown. 
Within the \emph{single parameter approximation} rationale, a unique scalar function depending on $\bmu$ is sought.
Following the strategy described for the spatial iteration, the high-order terms are neglected in the restriction of Equation~\eqref{eq:PGD-increment} to the parametric direction of the tangent manifold and $\De\phi$ is computed by solving the following algebraic equation
\begin{equation}\label{eq:PGD-param}
\left( \sum_{m=1}^n a_1^m \phi^m - a_2 \psi + a_3 \right) \De\phi 
= r_u^n + r_p^n ,
\end{equation}
where $r_u^n$ and $ r_p^n$ are the parametric residuals associated with the discretisation of the momentum and mass equations, respectively, and each coefficient $a_i, \, i{=}1,\ldots,3$ depends solely on the spatial functions $(\sigma_u^n\fu^n,\sigma_p^n\fp^n)$ and on the data of the problem, namely
\begin{equation}\label{eq:coef-param}
  \!\!\!\!\left\{\begin{aligned}
  a_1^m &:= \int_{V_i}{\!\!\! \sigma_u^n\fu^n {\cdot} \bigl[\grad {\cdot} (\sigma_u^n\fu^n {\otimes} \sigma_u^m\fu^m)\bigr] \, dV} \\
  				 & \hspace{20pt} + \int_{V_i}{\!\!\! \sigma_u^n\fu^n {\cdot} \bigl[\grad {\cdot} (\sigma_u^m\fu^m {\otimes} \sigma_u^n\fu^n)\bigr] \, dV} , \\
  a_2 &:=  \int_{V_i}{\!\!\! \sigma_u^n\fu^n {\cdot} \bigl[\grad {\cdot} (D \grad (\sigma_u^n\fu^n))\bigr] \, dV} , \\
  a_3 &:= \int_{V_i}{\!\!\! \sigma_u^n\fu^n {\cdot} \grad (\sigma_p^n\fp^n)  \, dV} 
			 + \!\! \int_{V_i}{\!\!\! \sigma_p^n\fp^n \grad {\cdot} (\sigma_u^n\fu^n) \, dV} .
  \end{aligned}\right.
\end{equation}
The unknown $\De\phi$ is discretised at the nodes of the parametric domain $\bI$ and the resulting algebraic equation is solved via a collocation method.
Similarly to the spatial iteration, the separable form of \eqref{eq:resU}-\eqref{eq:resP} is exploited to perform computationally efficient pointwise evaluations of the residuals at the nodes of $\bI$. 
The complete derivation of the separated form of the right-hand side is detailed in~\ref{sc:residuals}.

%==========================================================================
\subsection{A nonintrusive implementation of the proper generalised decomposition in OpenFOAM}
\label{sc:nonintrusive}
%==========================================================================

\begin{figure}[p]
\centering
	\includegraphics[width=1.0\textwidth]{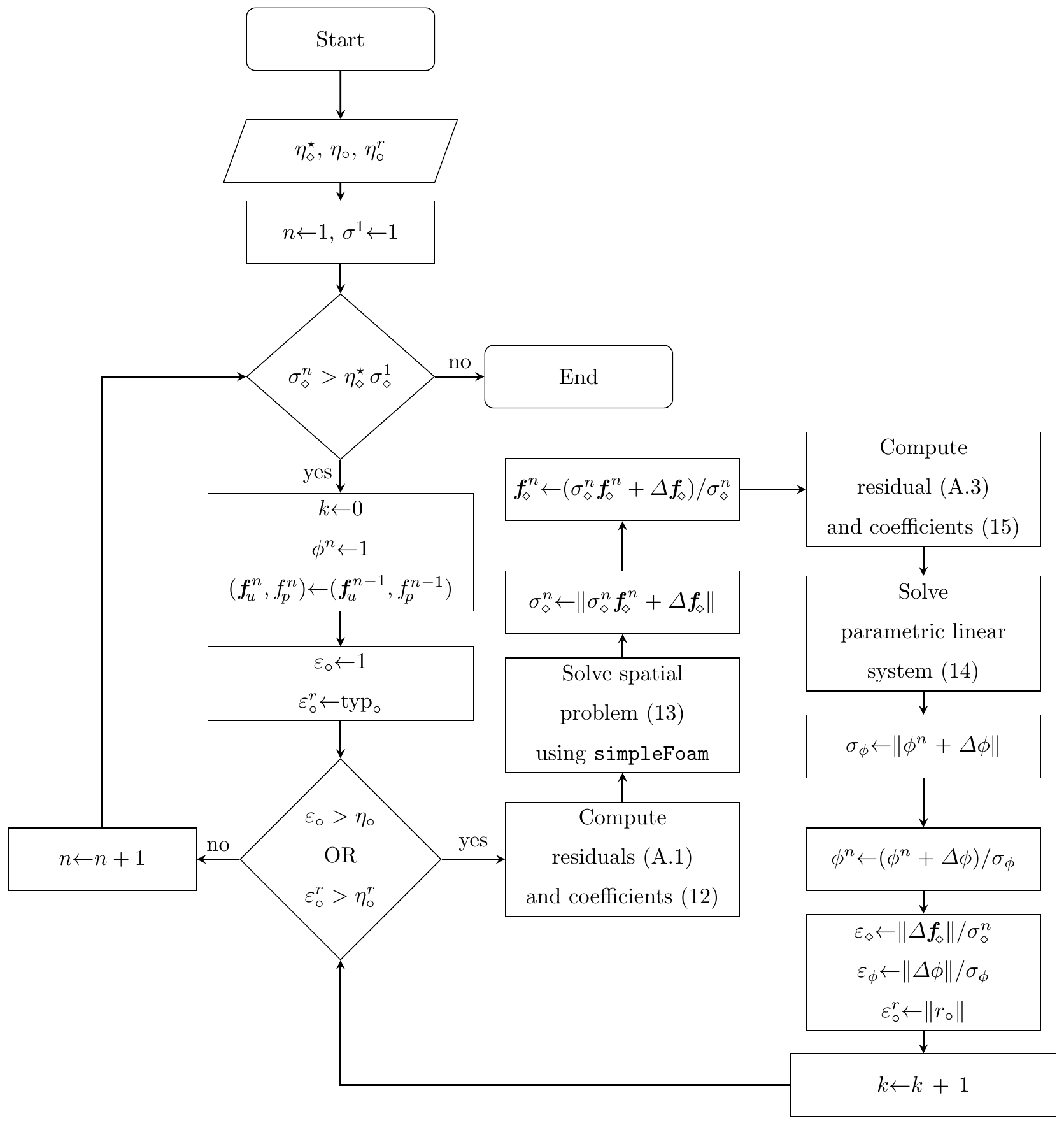}
	\caption{Flowchart of the nonintrusive \texttt{pgdFoam} algorithm. Legend: $\diamond=u,p$ and $\circ=u,p,\phi$.}
	\label{fig:flowchart}
\end{figure}

In order to construct an efficient PGD strategy applicable to engineering problems of interest for the industry, a critical aspect is its nonintrusiveness with respect to the OpenFOAM solving procedure \texttt{simpleFoam}.
As discussed in Section~\ref{sc:PGD-ADI}, inhomogeneous Dirichlet boundary conditions are treated by means of a spatial mode computed using the full-order solver, whereas the corresponding parametric mode is set equal to $1$ (Algorithm~\ref{alg:PGD-OF} - Step 1).
Then, the enrichment process is started and at each iteration of the alternating direction scheme a spatial mode is computed using \texttt{simpleFoam} (Algorithm~\ref{alg:PGD-OF} - Steps 7 to 10) and a linear system is solved to determine the corresponding parametric term (Algorithm~\ref{alg:PGD-OF} - Steps 11 to 14).
The alternating direction iterations stop when the computed corrections $\De\fd$, $\De\phi$ are negligible with respect to the amplitudes $\sigma_\diamond^n$, $\sigma_\phi$ of the current mode for $\diamond=u,p$ and the residuals $\varepsilon_\circ^r$ are sufficiently small for $\circ=u,p,\phi$ (Algorithm~\ref{alg:PGD-OF} - Steps 6 and 15).
The global enrichment strategy ends when the amplitude of the current mode $\sigma_\diamond^n$ is negligible with respect to the first one $\sigma_\diamond^1$ for $\diamond=u,p$ (Algorithm~\ref{alg:PGD-OF} - Step 3).
The complete flowchart of \texttt{pgdFoam} is displayed in Figure~\ref{fig:flowchart}.
\begin{remark}\label{rmrk:stop}
Alternative criterions may be considered to stop the greedy algorithm, e.g. when the magnitude of the last mode normalised with respect to the sum of the amplitudes of all the computed terms is lower than a user-defined tolerance $\eta_\diamond^\star$, namely
\begin{equation*}
\sigma_\diamond^n < \eta_\diamond^\star \sum_{m=1}^n\sigma_\diamond^m , \,\, \text{for} \,\, \diamond=u,p .
\end{equation*}
\end{remark}

\begin{algorithm}
\caption{\texttt{pgdFoam}: a nonintrusive PGD implementation in OpenFOAM}\label{alg:PGD-OF}
\begin{algorithmic}[1]
\REQUIRE{Tolerances $\eta_\diamond^\star$ for the greedy algorithm. Tolerances $\eta_\circ$ for the amplitudes and $\eta_\circ^r$ for the residuals in the alternating direction iteration. Typical values $\text{typ}_\circ$ for the residuals of the spatial and parametric problems. $\diamond=u,p$ and $\circ=u,p,\phi$.}
\STATE{Compute boundary condition modes: the spatial mode is solution of~\eqref{eq:weak-NS} using \texttt{simpleFoam} and the parametric mode is equal to $1$.}
\STATE{Set $n \gets 1$ and initialise the amplitudes of the spatial modes $\sigma_\diamond^1 \gets 1$.}

%___Loop modes (greedy)
\WHILE{$\sigma_\diamond^n > \eta_\diamond^\star\,\sigma_\diamond^1$}
\STATE{Set $k \gets 0$, the parametric predictor $\phi^n {\gets} 1$ and the spatial predictors $(\fu^n,\fp^n)$ using the last computed modes.}

\STATE{Initialise $\varepsilon_\circ \gets 1$, $\varepsilon_\circ^r \gets \text{typ}_\circ$.}
             
%___Loop ADI
\WHILE{$\varepsilon_\circ > \eta_\circ$ or $\varepsilon_\circ^r > \eta_\circ^r$}

\STATE{Compute the spatial residuals \eqref{eq:spatialRes} and coefficients \eqref{eq:coef-spatial}.}
\STATE{Solve the spatial Navier-Stokes problem \eqref{eq:PGD-spatial-final} using \texttt{simpleFoam}.}
\STATE{Normalise the spatial predictors: $\sigma_\diamond^n {\gets} \norm{\sigma_\diamond^n \fd^n + \De\fd}$.}
\STATE{Update the spatial predictors: $\fd^n {\gets} (\sigma_\diamond^n \fd^n + \De\fd)/\sigma_\diamond^n$.}
             
\STATE{Compute the parametric residual \eqref{eq:paramRes} and coefficients \eqref{eq:coef-param}.}
\STATE{Solve the parametric linear system \eqref{eq:PGD-param}.}
\STATE{Normalise the parametric predictor: $\sigma_\phi {\gets} \norm{\phi^n+\De\phi}$.}
\STATE{Update the parametric predictor: $\phi^n {\gets} (\phi^n + \De\phi)/\sigma_\phi$.}

\STATE{Update stopping criterions: $\varepsilon_\diamond {\gets} \norm{\De\fd}/\sigma_\diamond^n$, $\varepsilon_\phi {\gets} \norm{\De\phi}/\sigma_\phi$, $\varepsilon_\circ^r {\gets} \norm{r_\circ}$.}

\STATE{Update the alternating direction iteration counter: $k \gets k+1$.}

\ENDWHILE

\STATE{Update the mode counter: $n \gets n+1$.}
\ENDWHILE
\end{algorithmic}
\end{algorithm}

%==========================================================================
\section{Numerical validation}
\label{sc:simulations}
%==========================================================================

In this section, numerical examples are presented to validate the proposed methodology.
First, a test case with known analytical solution is considered to verify the optimal convergence rate of the high-dimensional FV approximation of the velocity and pressure fields, measured in the $\eltwo(\Omega {\times} \bI)$ norm, for a parametrised viscosity coefficient. 
In this context, special emphasis is given to the additional error introduced by the PGD, highlighting the range of applicability of the discussed reduced-order strategy in terms of expected accuracy of the parametric solution.
Moreover, a classical benchmark test for incompressible flow solvers, namely the nonleaky lid-driven cavity, is studied parametrising the imposed velocity of the lid in a range of values of the Reynolds number spanning from $1,\! 000$ to $4,\! 000$.

%==========================================================================
\subsection{Kovasznay flow with parametrised viscosity}
\label{sc:kovasznay}
%==========================================================================

Consider the Kovasznay flow~\cite{Kovasznay-Flow} for a parametrised viscosity $\nu(\mu){=}\mu$.
The analytical solution is
\begin{equation}\label{eq:KovasznayExact}
\begin{aligned}
\bu(x,y,\mu) &= \left(1-e^{\lambda(\mu) x}\cos(2\pi y) , \tfrac{\lambda(\mu)}{2\pi} e^{\lambda(\mu) x}\sin(2\pi y) \right) \\
p(x,y,\mu) &= \tfrac{1}{2} \left(1-e^{2\lambda(\mu) x} \right) + C
\end{aligned}
\end{equation}
where the constant $C$ is determined by fixing a reference value for the pressure field in one point of the domain, whereas $\lambda$ is a function of the parametrised viscosity and changes when the Reynolds number is modified, namely,
\begin{equation*}
\lambda(\mu) = \tfrac{1}{2\nu(\mu)} - \sqrt{ \tfrac{1}{(2\nu(\mu))^2} + 4\pi^2 } .
\end{equation*} 
The parameter $\mu$ is sought in the space $\I{=}[5 {\times} 10^{-3}, 10^{-2}]$, which is discretised with uniform intervals. The corresponding values of the Reynolds number span from $100$ to $200$.
The spatial domain $\Omega{=}[-1,1]^2$ is discretised with a family of Cartesian meshes of quadrilateral cells. 
The characteristic lengths $h_x$ and $h_{\mu}$ of the spatial and parametric discretisations, respectively, are provided in Table \ref{table:size1}.
\begin{table}[ht]
\centering
\begin{tabular}{|c || l | l | l | l | l | l |}
\hline
$h_x$    & $8.3\times 10^{-2}$ & $4{\times} 10^{-2}$ & $2{\times} 10^{-2}$ & $1{\times} 10^{-2}$     & $5{\times} 10^{-3}$      & $2.5{\times} 10^{-3}$ \\
\hline
$h_{\mu}$ & $2{\times} 10^{-2}$    & $1{\times} 10^{-2}$ & $5{\times} 10^{-3}$ & $2.5{\times} 10^{-3}$  & $1.25{\times} 10^{-3}$ & $6.25{\times} 10^{-4}$ \\
\hline
\end{tabular}
\caption{Normalised characteristic lengths of the spatial and parametric discretisations.}
\label{table:size1}
\end{table}

A convergence study under uniform mesh refinement is performed for the linearised Navier-Stokes equations using the meshes described in Table~\ref{table:size1}. 
In this context, a convective field $\ba$ given by the analytical expression of the Kovasznay velocity is introduced in Equation \eqref{eq:NavierStokes} and the convective term $\grad {\cdot} (\bu \otimes \bu)$ is replaced by $\grad {\cdot} (\bu \otimes \ba)$.
As detailed in Section \ref{sc:PGD-ADI}, an affine separation of the data is required to run PGD. 
Thus, the convective field $\ba$ is separated \emph{a priori} considering the first four terms of the Taylor expansion of $e^{\lambda x}$ in the analytical form of the velocity, see Equation~\eqref{eq:KovasznayExact}.
For $\mu{=}10^{-2}$, the relative $\mathcal{L}_2(\Omega)$ error of the resulting separated velocity field with respect to the exact one is $4.3 {\times} 10^{-3}$ and, consequently, a target error of $10^{-2}$ in the spatial discretisation is considered for the following convergence study.
Moreover, the Dirichlet boundary datum $\bu_D$ requires five modes to be described in a separated form. 

\begin{figure}%[ht]
\centering
	\includegraphics[width=0.6\textwidth]{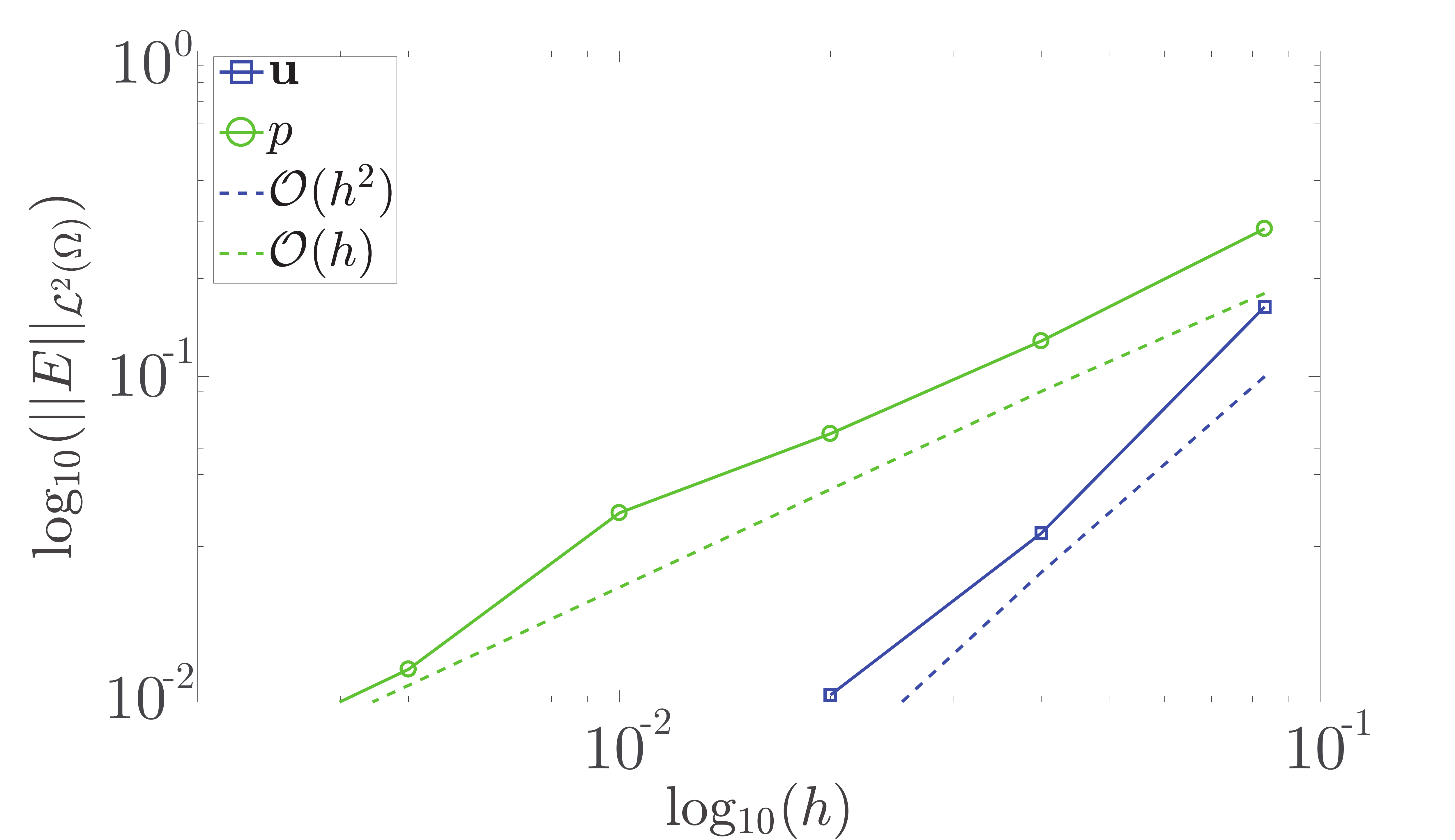}
	\caption{Optimal convergence of the $\eltwo(\Omega {\times} \bI)$ error of the PGD approximation of the Kovasznay flow with parametrised viscosity with respect to the exact solution as a function of the characteristic mesh size $h_x$.}
\label{fig:kovasznayError}
\end{figure}

The $\eltwo(\Omega {\times} \bI)$ error between the PGD approximation $(\upgd^n,\ppgd^n)$ computed using fifteen modes and the high-dimensional analytical solution $(\bu,p)$ as a function of the characteristic mesh size $h_x$ is displayed in Figure~\ref{fig:kovasznayError}. 
The optimal first-order convergence rate for pressure and second-order one for velocity are obtained.

To run the iterative procedure \texttt{pgdFoam} described in Algorithm~\ref{alg:PGD-OF}, a stopping criterion $\eta_{(u,p)}\leq 10^{-5}$ is considered, where $\eta_{(u,p)}$ accounts for the relative amplitude of both the velocity and pressure modes, namely
\begin{equation}\label{eq:relAmplitude}
\eta_{(u,p)} := \sqrt{\left(\frac{\sigma_u^n}{\sum_{m=1}^{n}\sigma_u^m}\right)^2+\left(\frac{\sigma_p^n}{\sum_{m=1}^{n}\sigma_p^m}\right)^2} .
\end{equation}
\begin{figure}%[ht]
\centering
	\subfigure[Amplitude of the computed modes]{\includegraphics[width=0.48\textwidth]{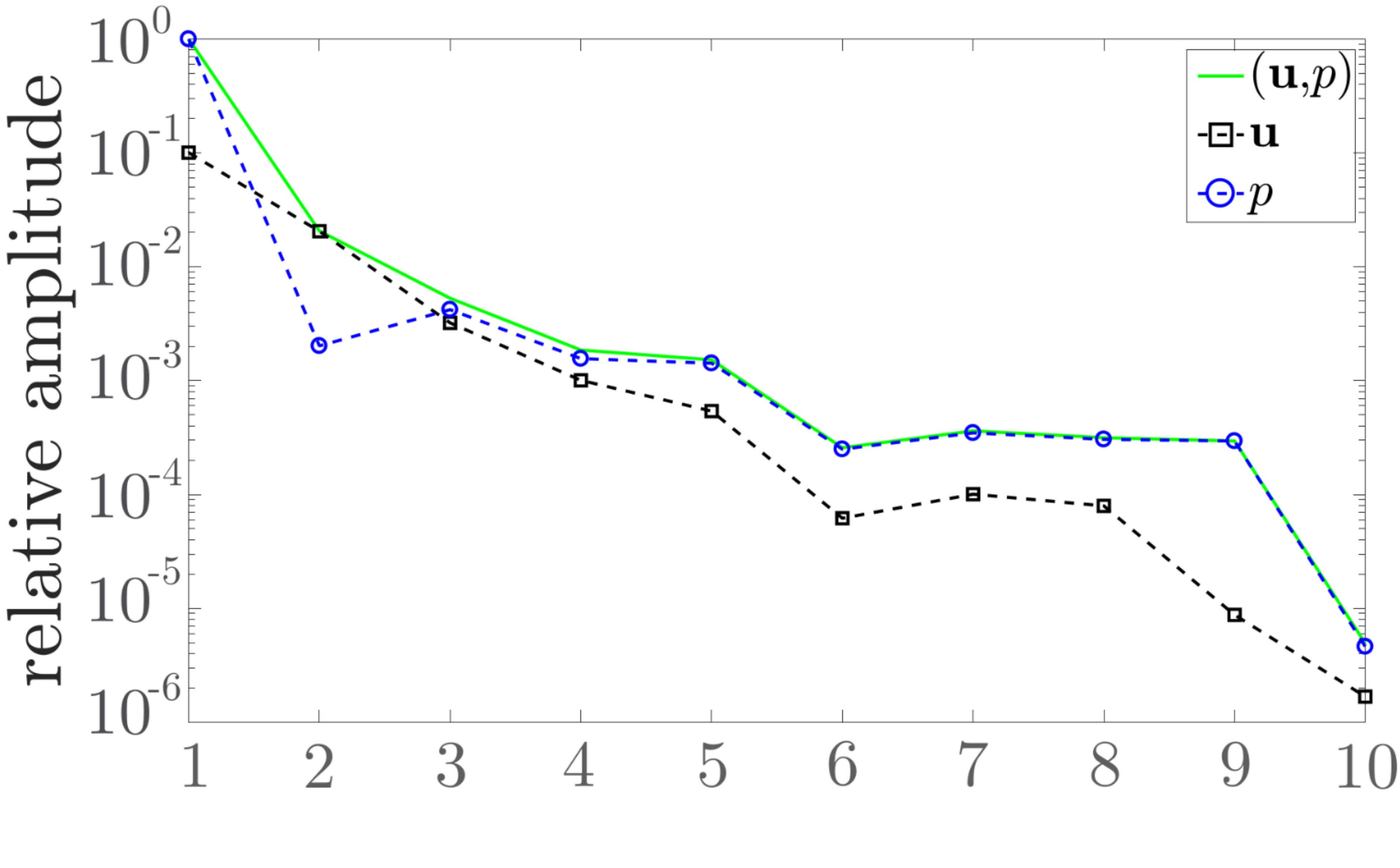}\label{fig:kovasznayAmplitude}}
	\subfigure[Computed parametric modes]{\includegraphics[width=0.48\textwidth]{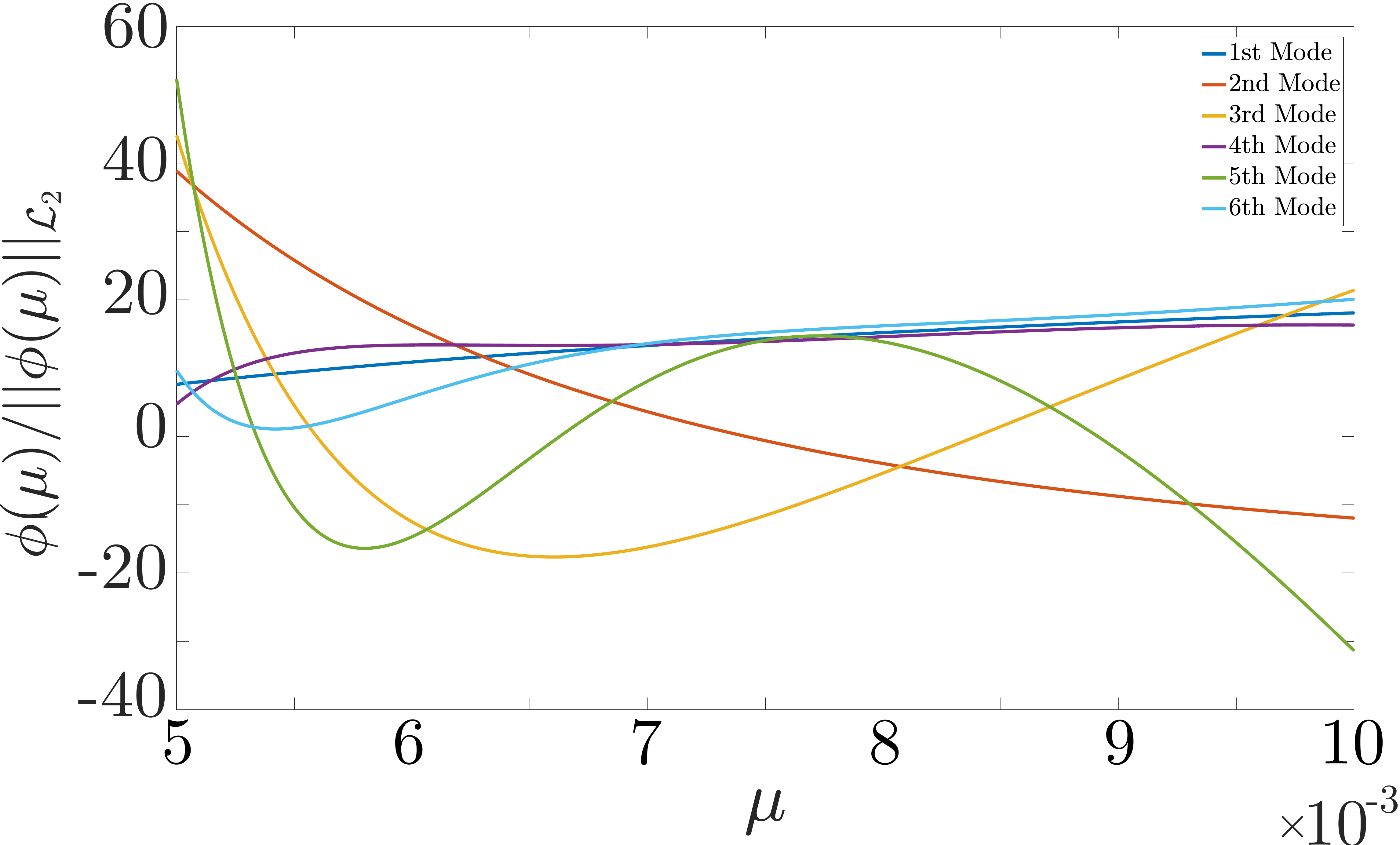}\label{fig:kovasznayParametricModes}}
	\caption{PGD approximation of the Kovasznay flow with parametrised viscosity. (a) Relative amplitude of the computed modes $\fu^m$ (black), $\fp^m$ (blue) and the combined amplitude of $(\fu^m,\fp^m)$ according to Equation~\eqref{eq:relAmplitude}. (b) First six normalised computed parametric modes.}
\label{fig:kovasznayAmplitude-ParametricModes}
\end{figure}
In Figure~\ref{fig:kovasznayAmplitude}, the evolution of the amplitude $\eta_{(u,p)}$, $\eta_u$ and $\eta_p$ is displayed for the finest mesh described in Table~\ref{table:size1}.
After ten computed modes, the stopping criterion is fulfilled and the PGD enrichment stops. 

\begin{figure}%[ht]
\centering
	\subfigure[$\fp^1$]{\includegraphics[width=0.15\textwidth]{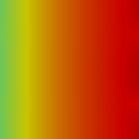}}
	\subfigure[$\fp^2$]{\includegraphics[width=0.15\textwidth]{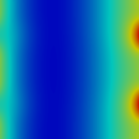}}
	\subfigure[$\fp^3$]{\includegraphics[width=0.15\textwidth]{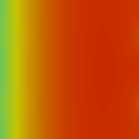}}
	\subfigure[$\fp^4$]{\includegraphics[width=0.15\textwidth]{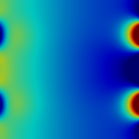}}
	\subfigure[$\fp^5$]{\includegraphics[width=0.15\textwidth]{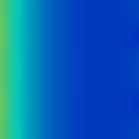}}
	\subfigure[$\fp^6$]{\includegraphics[width=0.15\textwidth]{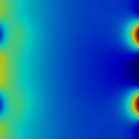}}
	
	\subfigure{\includegraphics[width=0.5\textwidth]{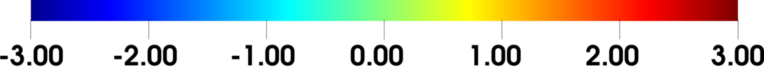}}
	\caption{PGD approximation of the Kovasznay flow with parametrised viscosity. First six computed spatial modes $\fp^m, \ m=1,\ldots,6$ for pressure.}
\label{fig:kovasznayPModes}
%\end{figure}
%
%\begin{figure}%[ht]
\centering
	\subfigure[$\fu^1$]{\includegraphics[width=0.15\textwidth]{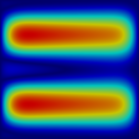}}
	\subfigure[$\fu^2$]{\includegraphics[width=0.15\textwidth]{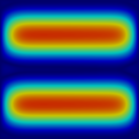}}
	\subfigure[$\fu^3$]{\includegraphics[width=0.15\textwidth]{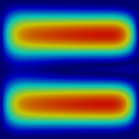}}
	\subfigure[$\fu^4$]{\includegraphics[width=0.15\textwidth]{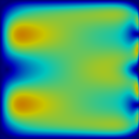}}
	\subfigure[$\fu^5$]{\includegraphics[width=0.15\textwidth]{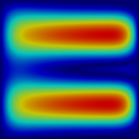}}
	\subfigure[$\fu^6$]{\includegraphics[width=0.15\textwidth]{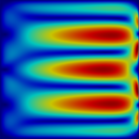}}

	\subfigure{\includegraphics[width=0.5\textwidth]{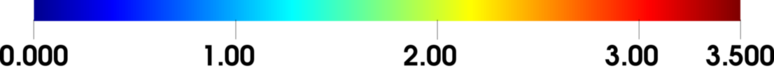}}
	\caption{PGD approximation of the Kovasznay flow with parametrised viscosity. First six computed spatial modes $\fu^m, \ m=1,\ldots,6$ for velocity.}
\label{fig:kovasznayUModes}
\end{figure}

As previously mentioned, five terms are required to describe the Dirichlet boundary conditions in a separated form. Henceforth, only the computed modes, starting from the sixth term of the PGD approximation are displayed.
In Figures~\ref{fig:kovasznayParametricModes}, the first six normalised computed parametric modes are displayed.
The corresponding computed spatial modes for pressure and velocity are presented in Figure~\ref{fig:kovasznayPModes} and \ref{fig:kovasznayUModes}, respectively.

\begin{figure}%[ht]
\centering
	\subfigure[$\ppgd^6$]{\includegraphics[width=0.2\textwidth]{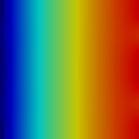}}
	\subfigure[$\ppgd^8$]{\includegraphics[width=0.2\textwidth]{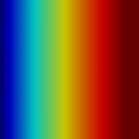}}	
	\subfigure[$\ppgd^{15}$]{\includegraphics[width=0.2\textwidth]{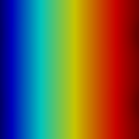}}
	\subfigure[Exact $p$]{\includegraphics[width=0.2\textwidth]{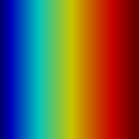}}
	\subfigure{\includegraphics[width=0.06\textwidth]{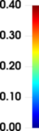}}
		
	\subfigure[$\upgd^6$]{\includegraphics[width=0.2\textwidth]{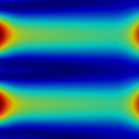}}
	\subfigure[$\upgd^8$]{\includegraphics[width=0.2\textwidth]{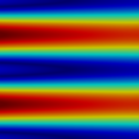}}
	\subfigure[$\upgd^{15}$]{\includegraphics[width=0.2\textwidth]{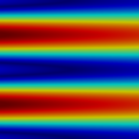}}
	\subfigure[Exact $\bu$]{\includegraphics[width=0.2\textwidth]{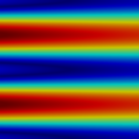}}
	\subfigure{\includegraphics[width=0.06\textwidth]{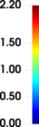}}
	\caption{Comparison of the PGD approximation to the analytical solution of the Kovasznay flow for $\text{Re}=200$, that is $\mu=10^{-2}$. Pressure (top) and velocity (bottom) using $6$, $8$ and $15$ terms, that is $1$, $3$ and $10$ computed modes besides the ones accounting for boundary conditions.}
\label{fig:kovasznayApprox}
\end{figure}

The PGD approximation $(\upgd^n,\ppgd^n)$ using $n{=}6,8,15$, that is with $1$, $3$ and $10$ computed modes, respectively, is compared to the analytical solution for the case of $\text{Re}{=}200$, in Figure~\ref{fig:kovasznayApprox}.

%==========================================================================
\subsection{Two-dimensional cavity with parametrised lid velocity}
\label{sc:cavity2D}
%==========================================================================

In this section, the classical benchmark problem of the nonleaky lid-driven cavity is studied~\cite{Ghia-Lid-Driven}.
The unitary square $\Omega{=}[0,1]^2$ is considered as spatial domain and homogeneous Dirichlet boundary conditions are imposed on the lateral and bottom walls. 
On the top wall, a velocity $\bu_{\text{lid}}(\bx,\mu){=}400\mu\bu_{\text{lid}}(\bx)$ is enforced, where the parameter $\mu\in[0.25, 1]$ acts as a scaling factor of the maximum velocity of the lid, whereas $\bu_{\text{lid}}(\bx)$ is a velocity profile featuring two ramps on the top-left and top-right corners of the domain to account for the change between null and maximum velocity.
As classical in the literature treating the lid-driven cavity example, for $x \in [0, 0.06]$ and $x \in [0.94, 1]$, the horizontal component of the lid velocity changes linearly from $0$ to $400\mu$ and vice versa.
The dynamic viscosity is set to $\nu{=}\SI{0.1}{m^2/s}$ and the values considered for the Reynolds number span from $1,\! 000$ to $4,\! 000$.

The nonlinear term of the Navier-Stokes equations is now treated as described in Section~\ref{sc:PGD-ADI}.
The mode handling the boundary conditions is obtained as a full-order solution of the Navier-Stokes equations using the \texttt{simpleFoam} algorithm for a lid velocity computed using $\mu{=}1$, that is for a maximum horizontal velocity of $\SI{400}{m/s}$.
The corresponding parametric boundary condition mode is set to be linearly evolving from $\mu{=}0.25$ to $\mu{=}1$, that is $\phi(\mu){=}\mu$.

\begin{figure}%[ht]
\centering
	\subfigure[Amplitude of the computed modes]{\includegraphics[width=0.46\textwidth]{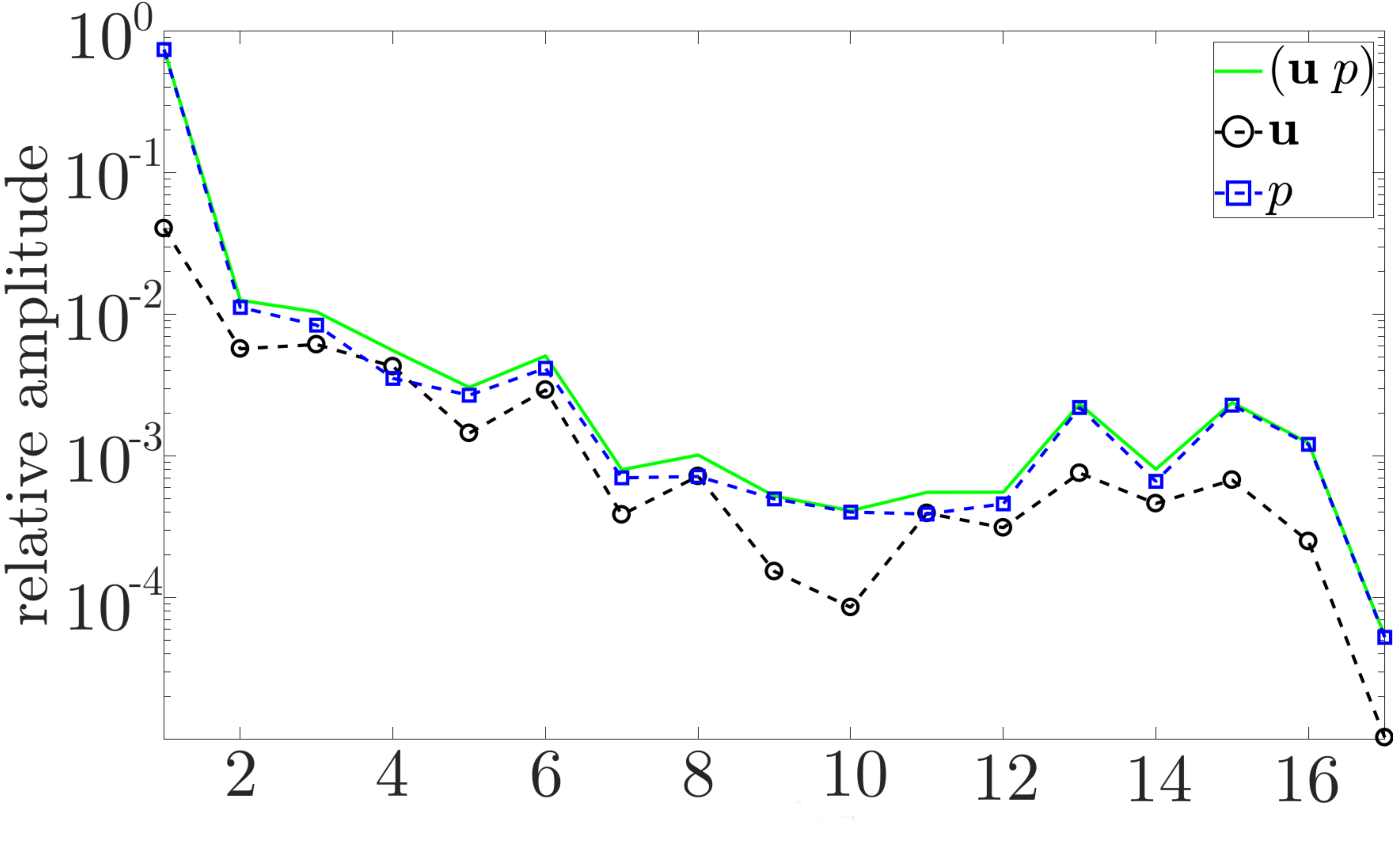}\label{fig:lidAmplitude}}
	\subfigure[Computed parametric modes]{\includegraphics[width=0.52\textwidth]{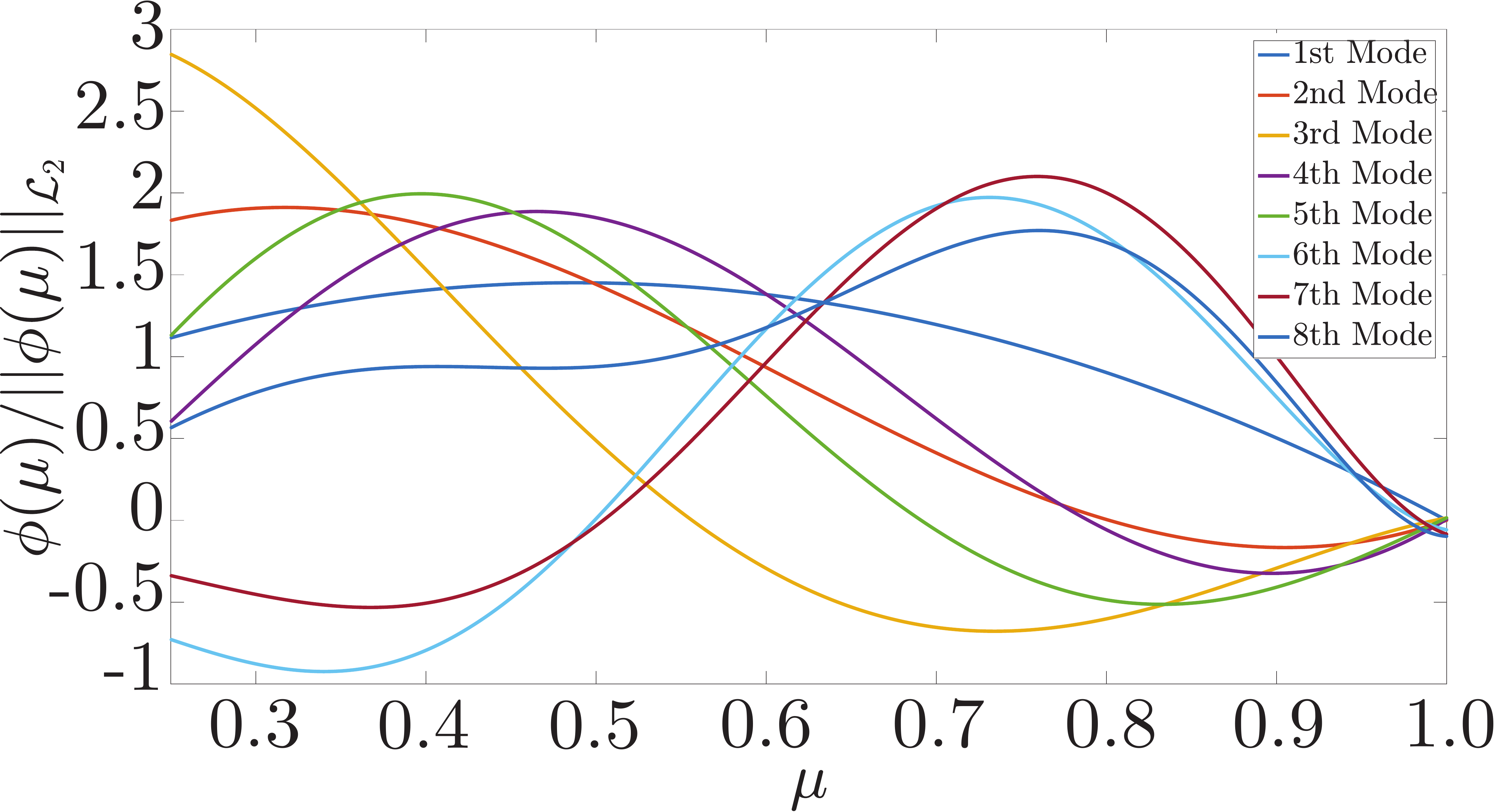}\label{fig:lidParametricModes}}
	\caption{PGD of the cavity flow with parametrised lid velocity. (a) Relative amplitude of the computed modes $\fu^m$ (black), $\fp^m$ (blue) and the combined amplitude of $(\fu^m,\fp^m)$ according to Equation~\eqref{eq:relAmplitude}. (b) First eight normalised computed parametric modes.}
\label{fig:lidAmplitude-ParametricModes}
\end{figure}

Following the rationale described in the previous section, two different stopping criterions are considered for the PGD enrichment strategy, namely $\eta_{(u,p)}\leq 10^{-3}$ and $\eta_{(u,p)}\leq 10^{-4}$.
Figure~\ref{fig:lidAmplitude} displays the relative amplitude of the computed modes. Note that the first stopping point is achieved after seven computed modes, whereas seventeen terms are required to fulfil the lower tolerance. 
The corresponding computed parametric modes are presented on Figure~\ref{fig:lidParametricModes}. It is worth noting that all the computed parametric modes are close or equal to $0$ for $\mu{=}1$. This is due to the fact that the boundary conditions of the problem are imposed by means of a full-order solution computed for the maximum value of $\mu$ in the parametric space. Hence, the case of $\mu{=}1$ is accurately described by the PGD approximation using solely the mode obtained via \texttt{simpleFoam}.

\begin{figure}%[ht]
\centering
	\includegraphics[width=0.6\textwidth]{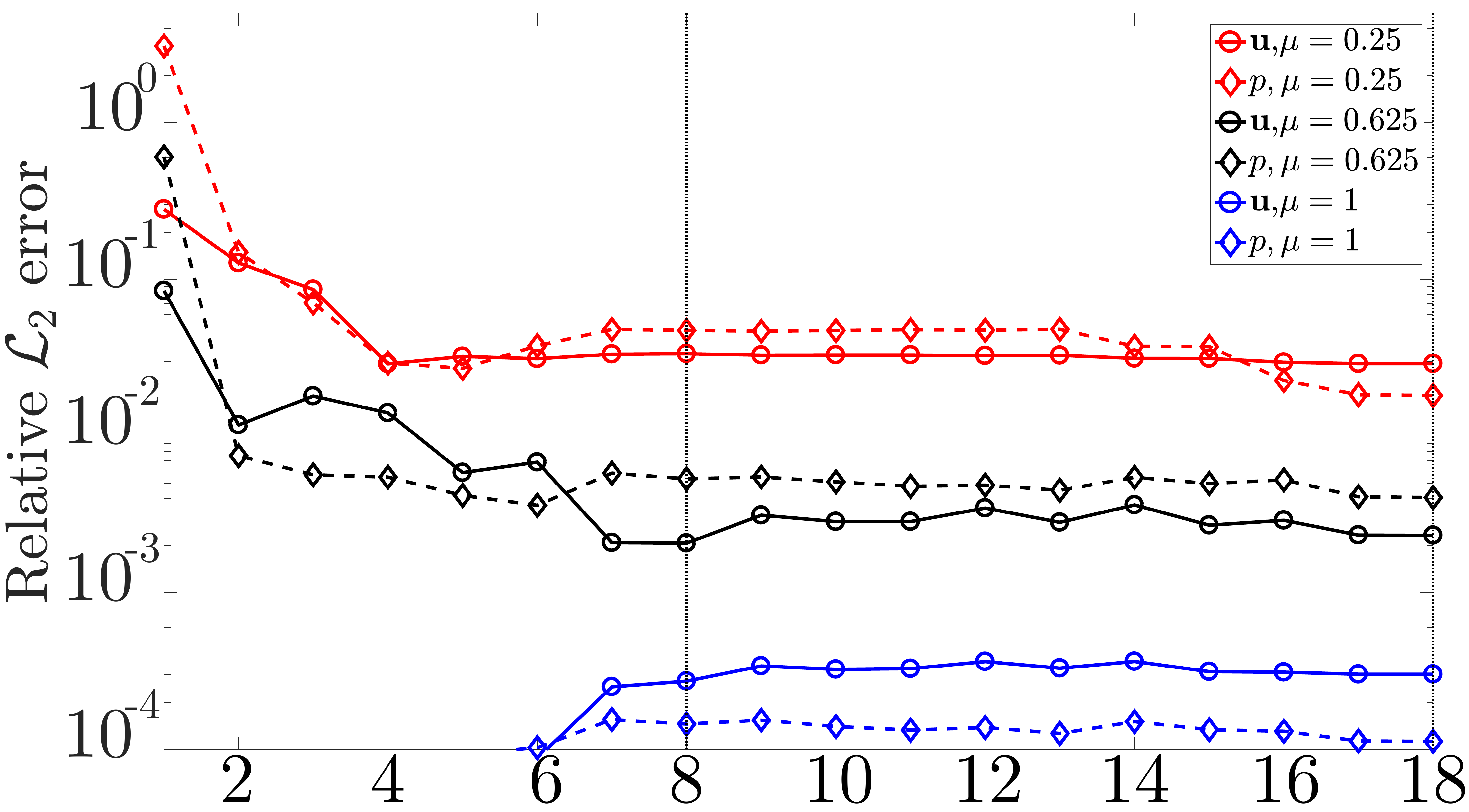}
	\caption{Relative $\mathcal{L}_2(\Omega)$ errors of the PGD approximation of the cavity flow with parametrised lid velocity with respect to the full-order solution as a function of the global number of modes (i.e. boundary conditions and computed) utilised in the PGD expansion.}
\label{fig:lid_errors}
\end{figure}

Now, the online evaluations of the PGD approximation of the velocity and pressure fields for different values of the parameter $\mu$ are compared to the full-order solutions computed using \texttt{simpleFoam}. The corresponding relative $\eltwo(\Omega)$ errors are presented in Figure~\ref{fig:lid_errors} as a function of the number of modes utilised in the PGD approximation.
The first three computed modes, for which $\eta_{(u,p)}\leq 10^{-2}$, provide a good approximation of both velocity and pressure and limited corrections are introduced by the following modes until the stopping criterion of $10^{-3}$ is fulfilled at the vertical dotted line.
For the case $\mu{=}1$, a small error of the order of $10^{-4}$ appears starting from the fifth computed mode, i.e. $n{=}6$. This is due to the fact that the boundary condition mode already captures all the features of the flow, being a full-order solution of the Navier-Stokes equations as previously mentioned.

\begin{figure}%[ht]
	\centering
	\begin{tabular}[c]{@{}c@{}c@{ }c@{ }c@{ }}
		$\upgd$ & 
		\parbox[c]{0.3\textwidth}{\includegraphics[width=0.3\textwidth]{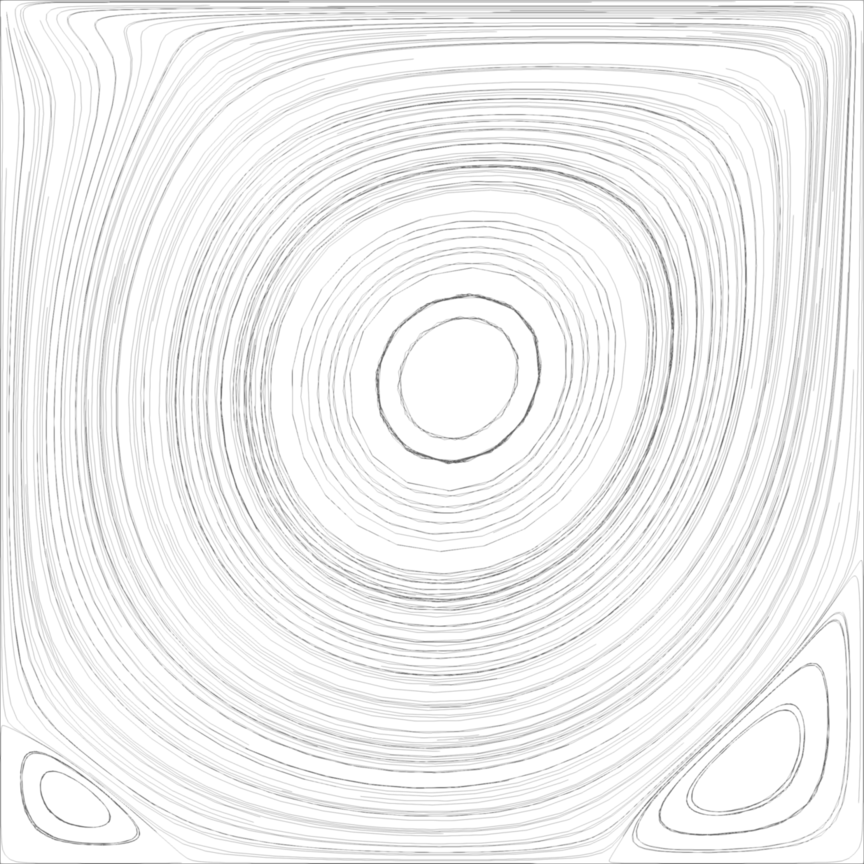}} & 
	 	\parbox[c]{0.3\textwidth}{\includegraphics[width=0.3\textwidth]{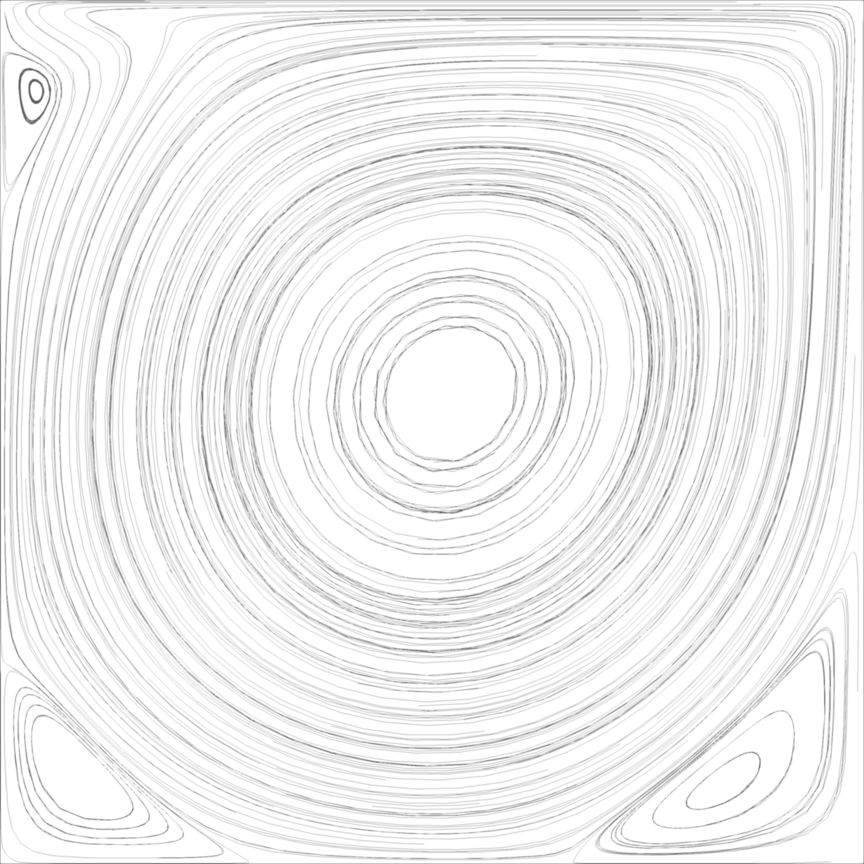}} &
		\parbox[c]{0.3\textwidth}{\includegraphics[width=0.3\textwidth]{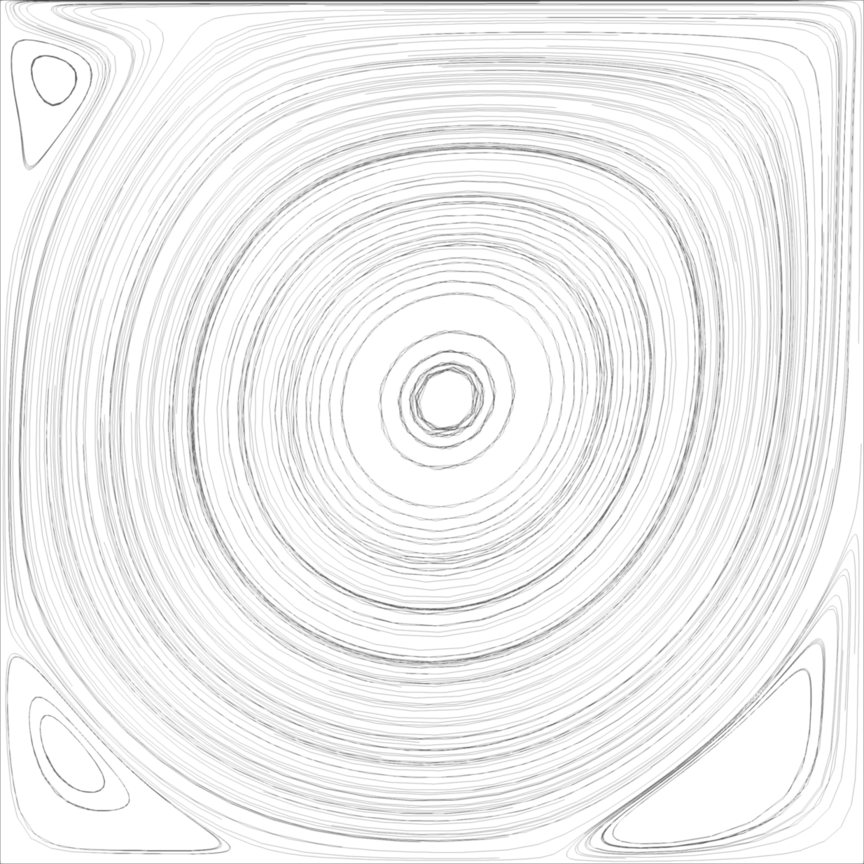}} \\[2em]
		$\uref$ &
		\parbox[c]{0.3\textwidth}{\includegraphics[width=0.3\textwidth]{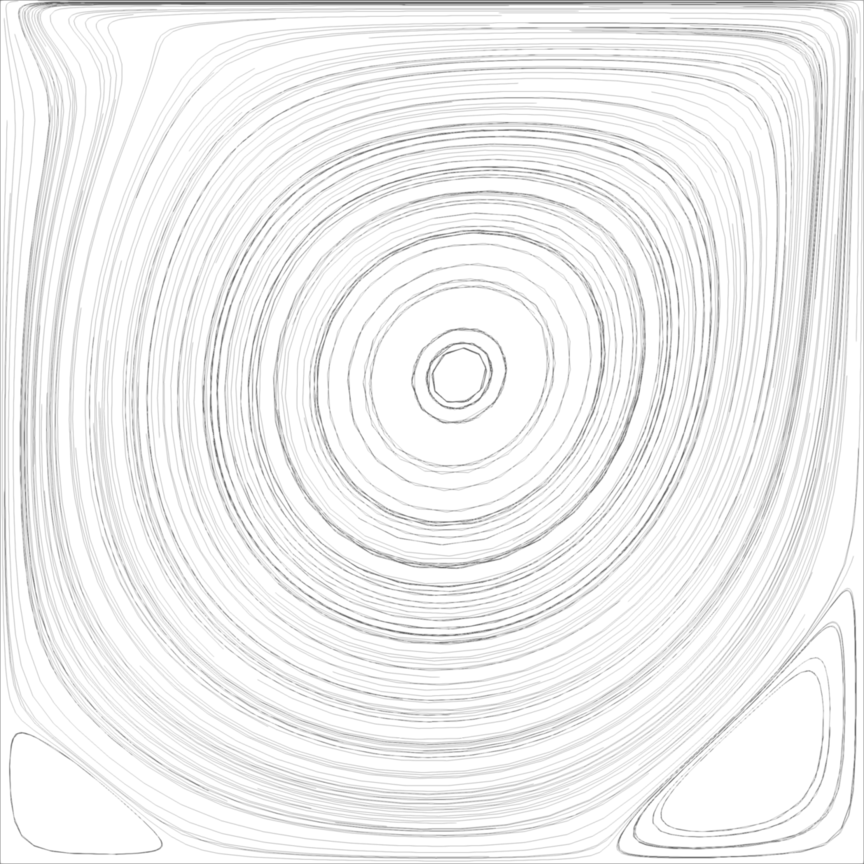}} &
		\parbox[c]{0.3\textwidth}{\includegraphics[width=0.3\textwidth]{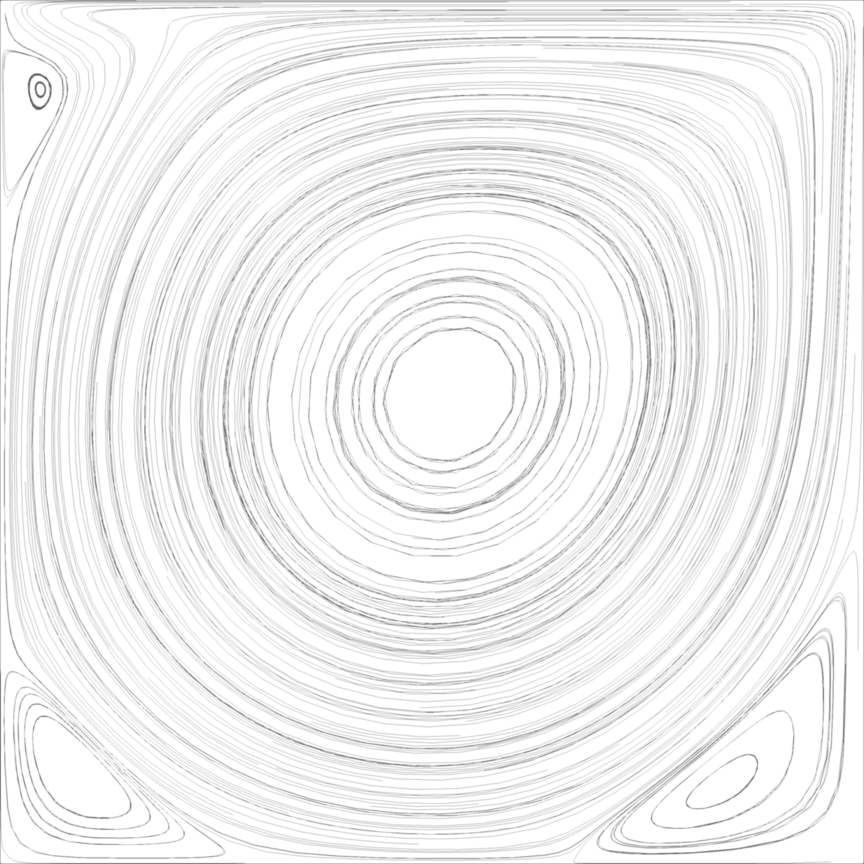}} &
		\parbox[c]{0.3\textwidth}{\includegraphics[width=0.3\textwidth]{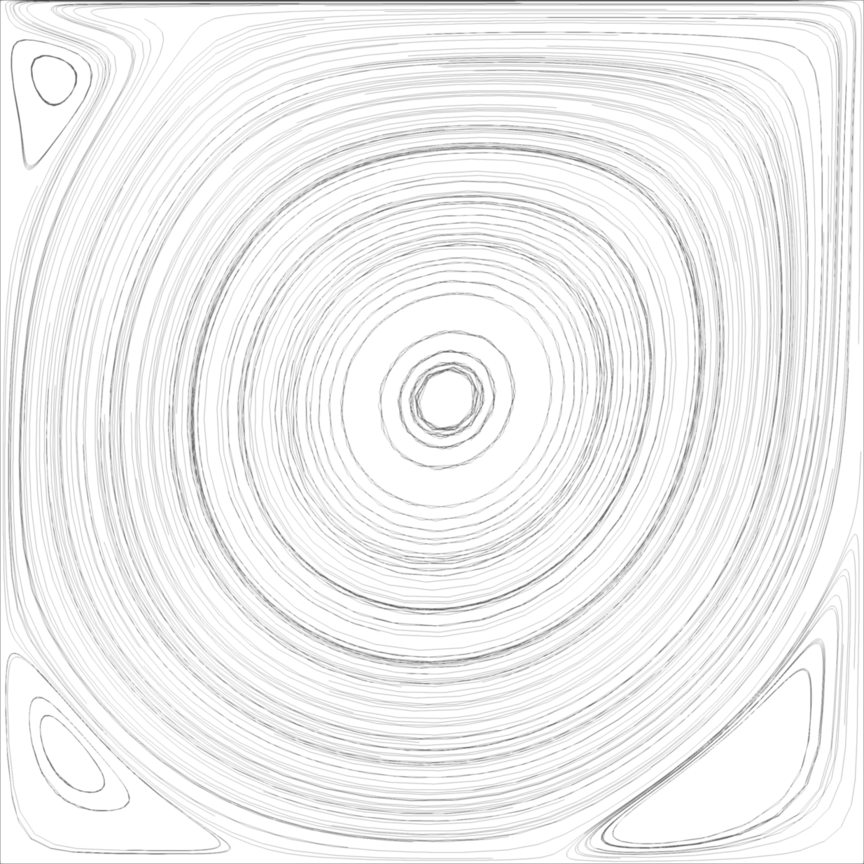}} \\
		  & $\mu{=}0.25$ & $\mu{=}0.625$ & $\mu{=}1$ 
	\end{tabular}
	\caption{Comparison of the PGD approximation (top) and the full-order solution (bottom) of the parametrised lid-driven cavity problem for maximum velocity of the lid of $100$, $250$ and $\SI{400}{m/s}$.}
\label{fig:lid_streamlines}
\end{figure}

A qualitative comparison of the reduced-order and full-order solutions of the parametrised lid-driven cavity problem is displayed in Figure~\ref{fig:lid_streamlines}.
Using the first seven computed modes, the PGD approximations for $\mu{=}0.25$, $\mu{=}0.625$ and $\mu{=}1$ are presented as long as their corresponding simulations obtained using \texttt{simpleFoam}. The cases under analysis correspond to a maximum horizontal velocity of the lid of $100$, $250$ and $\SI{400}{m/s}$, respectively.
It is worth noting that \texttt{pgdFoam} is able to capture the topological changes of the flow with great accuracy, managing to identify location and size of the vortices, as well as their appearance and disappearance according to the values of the Reynolds number considered in the analysis.

%==========================================================================
\section{Application to parametrised flow control problems}
\label{sc:flowControl}
%==========================================================================

Dynamically controlling the features of a flow is a challenging problem with several high-impact applications including, e.g., drag minimisation, stall control and aerodynamic noise reduction~\cite{Duvigneau-DV:06,Dede:07,Duvigneau-GDL:14}.
A major bottleneck to the design of flow control devices is represented by the large number of simulations involved in the tuning of the control loop.
In this section, the potential of the described nonintrusive PGD implementation in OpenFOAM is demonstrated for parametrised flow control problems.
Two- and three-dimensional internal flows with blowing jets are studied.
Specifically, a parametric study involving the peak velocity of the jets as extra-coordinate of the problem is considered to test the proposed PGD methodology.

%==========================================================================
\subsection{Lid-driven cavity with parametrised jet velocity}
\label{sc:jet2D}
%==========================================================================

Consider the nonleaky lid-driven cavity problem introduced in Section~\ref{sc:cavity2D}. 
The lid velocity is defined with two linear ramps, increasing from $0$ to $\SI{10}{m/s}$ on the top-left corner and decreasing correspondingly on the top-right one.
Three jets of size $\SI{0.12}{m}$ are introduced on the vertical walls, two on the right wall and one on the left, respectively. 
The parametrised velocity of the jets is $\bu_{\text{jet}}(\bx,\mu){=}\mu\bu_{\text{jet}}(\bx)$, where the maximum velocity is controlled by the parameter $\mu\in[0,1]$ and the profile $\bu_{\text{jet}}(\bx)$ is defined as
\begin{equation}\label{eq:lid_jets_BCs}
\bu_{\text{jet}}(x,y){=}\!
\begin{cases}\!
\Bigl({-}1{-}\cos\bigl(2\pi y/0.12\bigr), 0\Bigr)             &\!\!\text{for $x{=}1$, $y \in [0, 0.12]$,} \\ \!
\Bigl({-}1{+}\cos\bigl({-}2\pi(y{-}0.88)/0.12\bigr), 0\Bigr) &\!\!\text{for $x{=}1$, $y \in [0.88, 1]$,} \\ \!
\Bigl({-}1{+}\cos\bigl({-}2\pi(y{-}0.88)/0.12\bigr), 0\Bigr) &\!\!\text{for $x{=}0$, $y \in [0.88, 1]$.}
\end{cases} 
\end{equation}
An outlet boundary is added on the left vertical wall for $y \in [0, 0.12]$ and a free-traction condition is enforced.
The dynamic viscosity is set to $\nu{=}\SI{0.01}{m^2/s}$, therefore the corresponding Reynolds number is $\text{Re}{=}1,\! 000$. 

\begin{figure}%[ht]
	\centering
%	\subfigure[B.C. mode for the lid]{\includegraphics[width=0.351\textwidth]{lidJets-uBC0}}
	\subfigure[B.C. mode for the lid]{\includegraphics[width=0.42\textwidth]{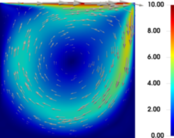}}
%	\subfigure[B.C. mode for the jets]{\includegraphics[width=0.3978\textwidth]{lidJets-uBC1}}
	\subfigure[B.C. mode for the jets]{\includegraphics[width=0.45\textwidth]{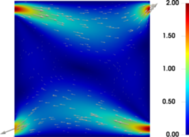}}
	\caption{Cavity flow with parametrised jet velocity. Spatial boundary condition modes for velocity.}
\label{fig:lidJets_BCmodes}
\end{figure}

The boundary conditions of the problem are enforced through two modes computed as full-order solutions via \texttt{simpleFoam} as shown in Figure~\ref{fig:lidJets_BCmodes}: the first one, for $\mu{=}0$, corresponds to lid velocity of $\SI{10}{m/s}$ and inactive jets; the second one, for $\mu{=}1$, is associated with the maximum velocity of the jets and a zero velocity of the lid. The corresponding parametric modes for the boundary conditions are set to $\phi(\mu){=}1$ and $\phi(\mu){=}\mu$, respectively.

\begin{figure}%[ht]
	\centering
	\begin{tabular}[c]{@{}c@{}c@{ }c@{ }c@{ }c@{ }c@{ }c@{}}
		$\upgd$ & 
		\parbox[c]{0.15\textwidth}{\includegraphics[width=0.15\textwidth]{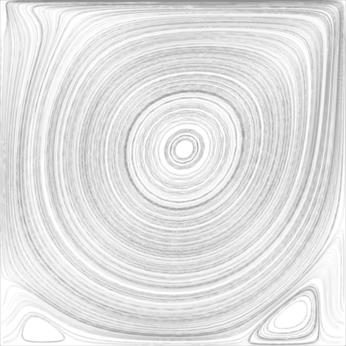}} & 
	 	\parbox[c]{0.15\textwidth}{\includegraphics[width=0.15\textwidth]{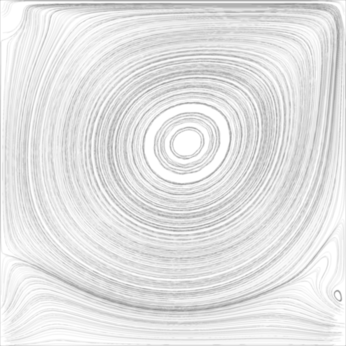}} &
		\parbox[c]{0.15\textwidth}{\includegraphics[width=0.15\textwidth]{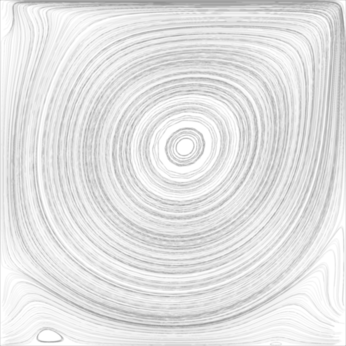}} &
		\parbox[c]{0.15\textwidth}{\includegraphics[width=0.15\textwidth]{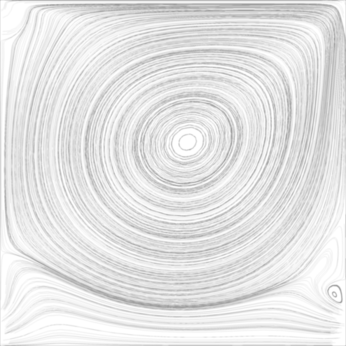}} &
		\parbox[c]{0.15\textwidth}{\includegraphics[width=0.15\textwidth]{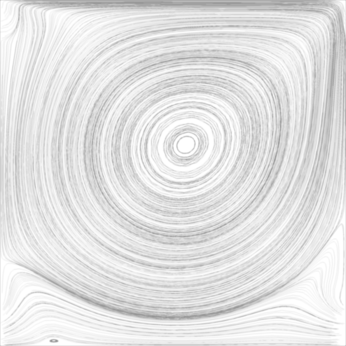}} &
		\parbox[c]{0.15\textwidth}{\includegraphics[width=0.15\textwidth]{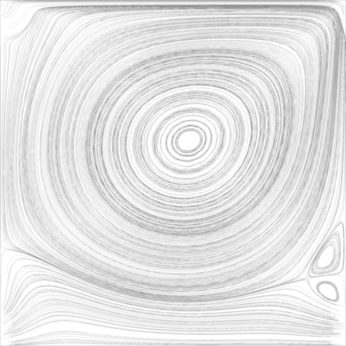}} \\[2em]
		$\uref$ &
		\parbox[c]{0.15\textwidth}{\includegraphics[width=0.15\textwidth]{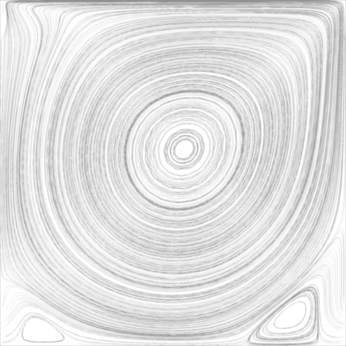}} &
		\parbox[c]{0.15\textwidth}{\includegraphics[width=0.15\textwidth]{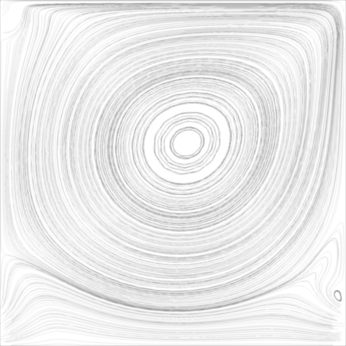}} &
		\parbox[c]{0.15\textwidth}{\includegraphics[width=0.15\textwidth]{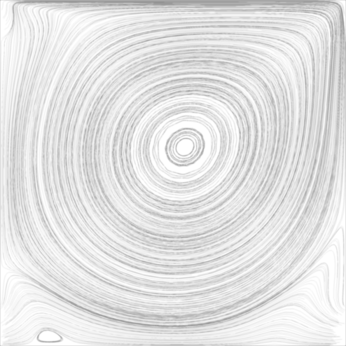}} &
		\parbox[c]{0.15\textwidth}{\includegraphics[width=0.15\textwidth]{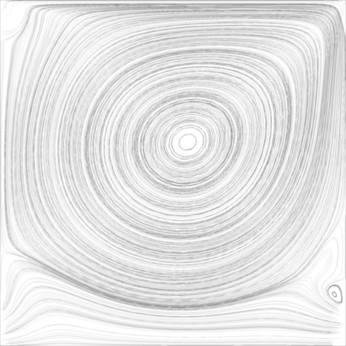}} &
		\parbox[c]{0.15\textwidth}{\includegraphics[width=0.15\textwidth]{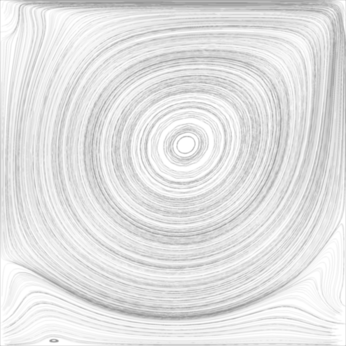}} &
		\parbox[c]{0.15\textwidth}{\includegraphics[width=0.15\textwidth]{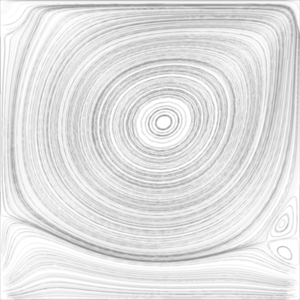}} \\
		  & $\mu{=}0.1$ & $\mu{=}0.3$ & $\mu{=}0.5$ & $\mu{=}0.7$ & $\mu{=}0.8$ & $\mu{=}1$\\
	\end{tabular}
\caption{Comparison of the PGD approximation (top) and the full-order solution (bottom) of the lid-driven cavity with jets for $\mu{=}0.1$, $\mu{=}0.3$, $\mu{=}0.5$, $\mu{=}0.7$, $\mu{=}0.8$ and $\mu{=}1$. }
\label{fig:lidJets_streamlines}
\end{figure}

Following the rationale previously discussed, the PGD enrichment process is stopped when $\eta_{(u,p)} \leq 10^{-4}$.
%
%In Figure~\ref{fig:lidJets_amplitude}, the amplitude of the computed modes is displayed.
%%
%\begin{figure}[ht]
%	\centering
%	\includegraphics[width=0.5\textwidth]{LidJets-amplitudeLidJets}
%	
%	\caption{Relative amplitude of the computed modes $\fu^m$ (black), $\fp^m$ (blue) and the combined amplitude of $(\fu^m,\fp^m)$ according to Equation~\eqref{eq:relAmplitude}.}
%	\label{fig:lidJets_amplitude}
%\end{figure}
%
In Figure~\ref{fig:lidJets_streamlines}, the generalised solution computed by the PGD is particularised for several values of the parameter under analysis and compared with the corresponding full-order solutions provided by \texttt{simpleFoam}.
The flows for $\mu{=}0.1$, $\mu{=}0.3$, $\mu{=}0.5$, $\mu{=}0.7$, $\mu{=}0.8$ and $\mu{=}1$ are displayed, covering a wide range of flow regimes in the cavity.
It is worth noting that the discussed reduced-order strategy is able to capture the topological changes in the flow features and accurately reproduce the appearance and disappearance of vortices in different regions of the domain.

\begin{figure}%[ht]
\centering
	\includegraphics[width=0.7\textwidth]{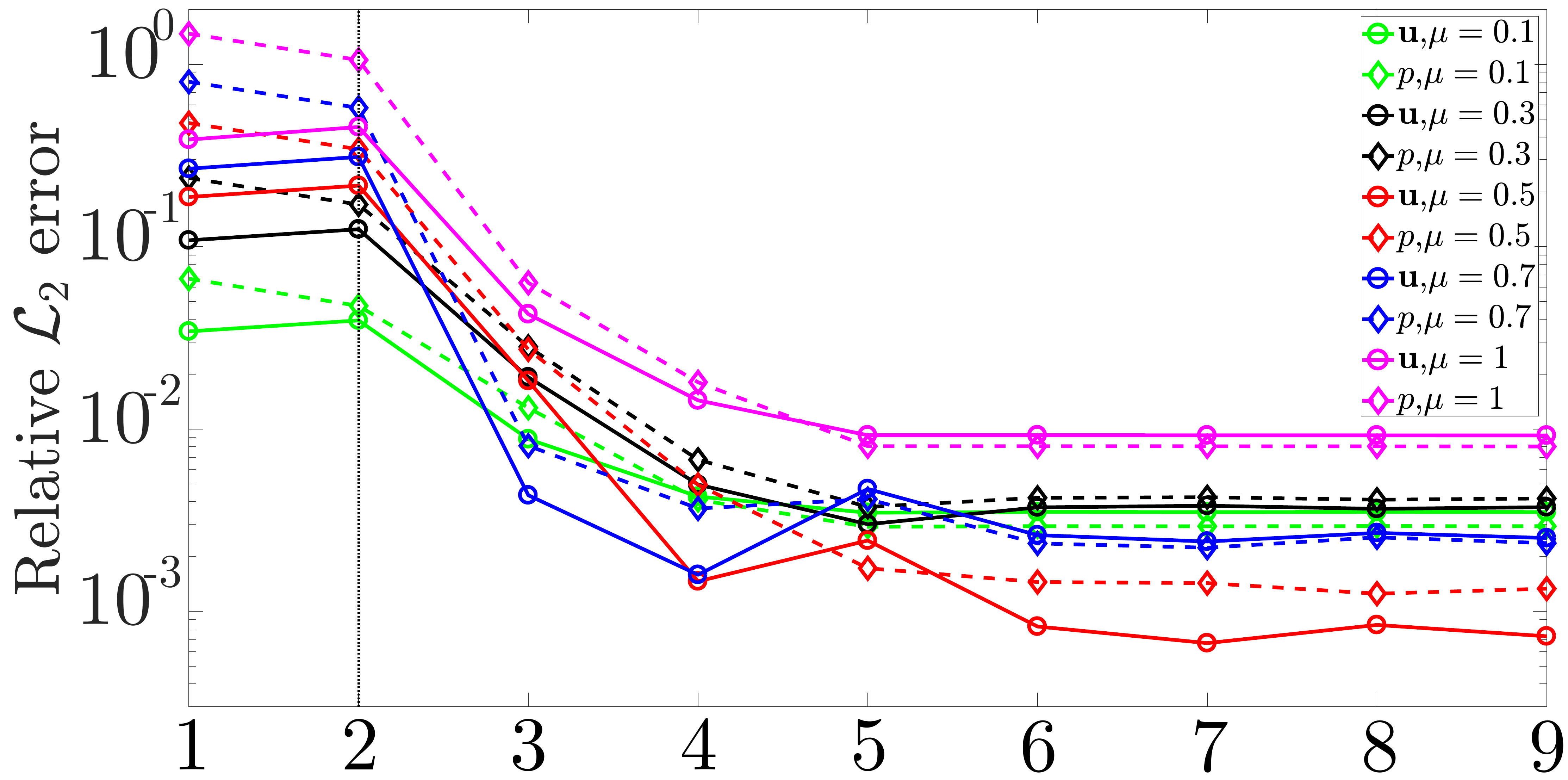}
	\caption{Relative $\eltwo(\Omega)$ error of the PGD approximation of the cavity flow with parametrised jet velocity with respect to the full-order solution for $\mu{=}0.1$, $\mu{=}0.5$, $\mu{=}0.7$ and $\mu{=}1$. The vertical dotted line separates the first two modes accounting for the boundary conditions and the computed modes.}
\label{fig:lidJets_errors}
\end{figure}

The accuracy of the PGD approximation with respect to the full-order solution is also verified by computing the relative $\eltwo(\Omega)$ error of the spatial discretisation while enriching the modal description of the solution. Specifically, Figure~\ref{fig:lidJets_errors} shows that using seven computed modes all approximations present relative errors lower than $10^{-2}$. 

%==========================================================================
\subsection{S-Bend with flow control driven by a jet}
\label{sc:duct}
%==========================================================================

\begin{figure}%[ht]
	\centering
	\subfigure[Front view]{\includegraphics[width=0.4\textwidth]{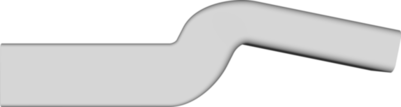}}
	\subfigure[Bottom view]{\includegraphics[width=0.4\textwidth]{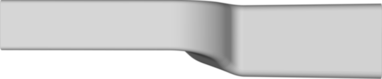}}
	
	\subfigure[Perspective view]{\includegraphics[width=0.4\textwidth]{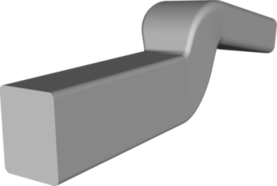}}
	\subfigure[Patches of the duct]{\includegraphics[width=0.4\textwidth]{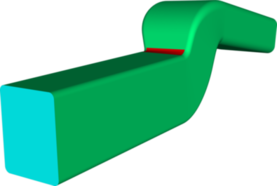}}

\caption{Geometrical model of the S-Bend. On the bottom-right image, the jet patch is highlighted in red.}
\label{fig:sBendGeom}
\end{figure}

In this section, the proposed PGD methodology is applied to a flow control problem using a three-dimensional geometry of industrial interest.
The model of a heating, ventilation and air conditioning (HVAC) duct section provided by Volkswagen AG is shown in Figure~\ref{fig:sBendGeom}.
A jet is introduced on the red patch, at the first bend of the duct.
The velocity profile of the jet is a sinusoidal function defined on the reference planar square $[0,1]^2$ as
\begin{equation}\label{eq:sBend_jet_BCs}
\begin{aligned}
&u_{\hat{y}}(\hat{x},\hat{z})= 0.0375(1-\cos(-2\pi\hat{x}))(1-\cos(2\pi\hat{z})) 
\end{aligned}
\end{equation}
and pointing in the direction $\hat{y}$ orthogonal to the plane $(\hat{x},\hat{z})$.
The parametrisation is constructed as a scaling of the jet velocity from $u_y{=}\SI{-0.015}{m/s}$, i.e. blowing, to suction with $u_y{=}\SI{0.15}{m/s}$.
A single parameter $\mu$ is introduced and the parametric domain considered for the analysis is $\I{=}[-0.1,1]$. Note that this problem is particularly challenging due to the change of sign in the interval of parametric values considered leading to different physical phenomenons.
The remaining boundary conditions feature homogeneous velocity on all the lateral walls, a parabolic velocity profile with mean value $\bu{=}(0.83,0,0)\SI{}{m/s}$ on the inlet and a free-traction on the outlet.
The dynamic viscosity is set to $\nu{=}1.588{\times} 10^{-5} \, \SI{}{m^2/s}$ and the corresponding value of the Reynolds number is $\text{Re}{=}280$. 
The quantity of interest in this problem is the pressure drop computed along the duct. 

As previously done for the lid-driven cavity with jets, two modes to account for the boundary conditions are computed using \texttt{simpleFoam}. The first mode is a full-order solution corresponding to the case of inactive jet and given inlet parabolic profile; the second one, is obtained setting a zero inlet velocity and a jet of maximum velocity $u_y{=}\SI{0.15}{m/s}$. The corresponding parametric modes are $\phi(\mu){=}1$ and $\phi(\mu){=}\mu$, respectively.

\begin{figure}%[ht]
\centering
	\subfigure[Amplitude of the computed modes]{\includegraphics[width=0.48\textwidth]{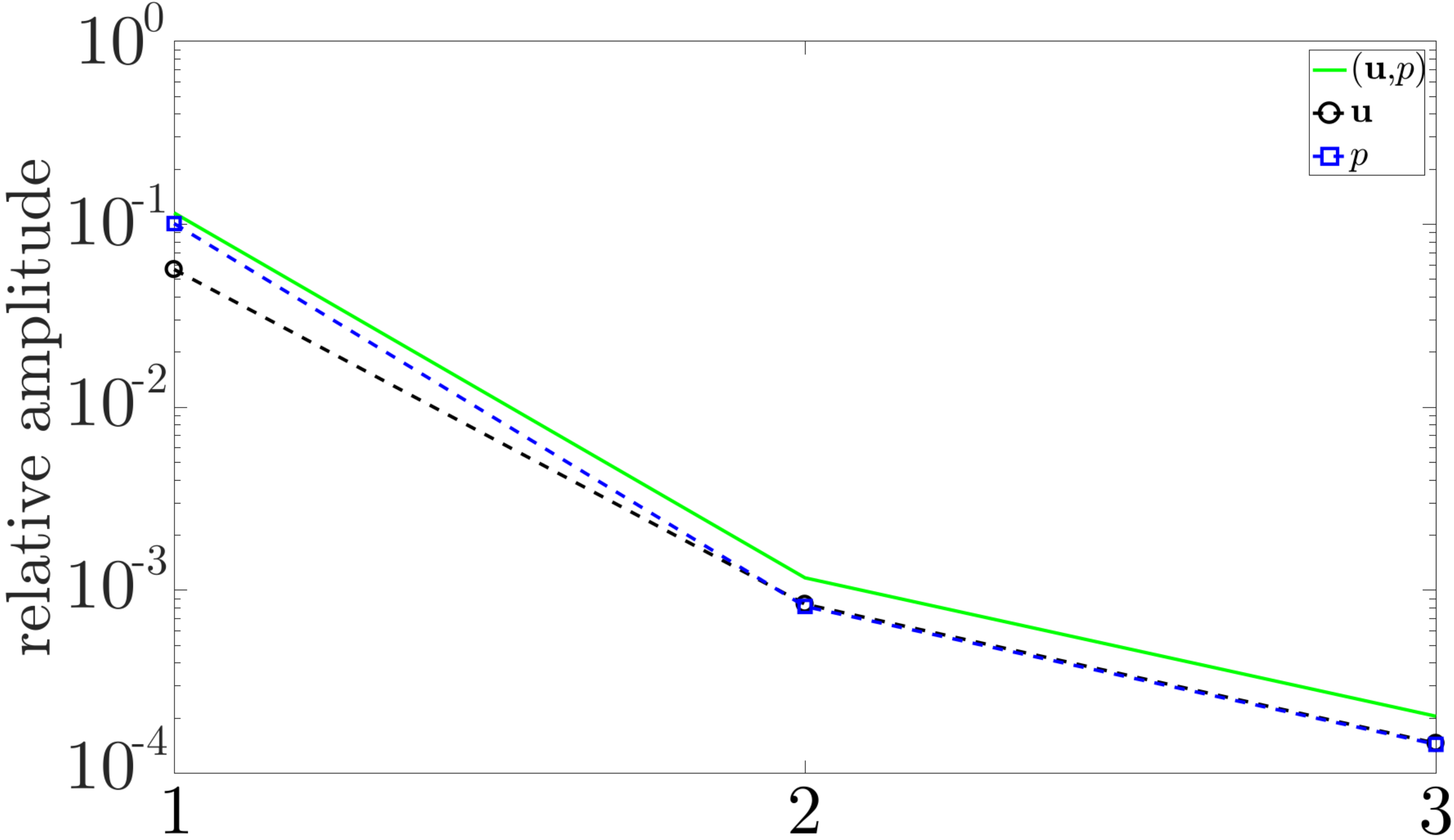}\label{fig:sBendAmplitude}}
	\subfigure[Relative $\eltwo(\Omega)$ error]{\includegraphics[width=0.48\textwidth]{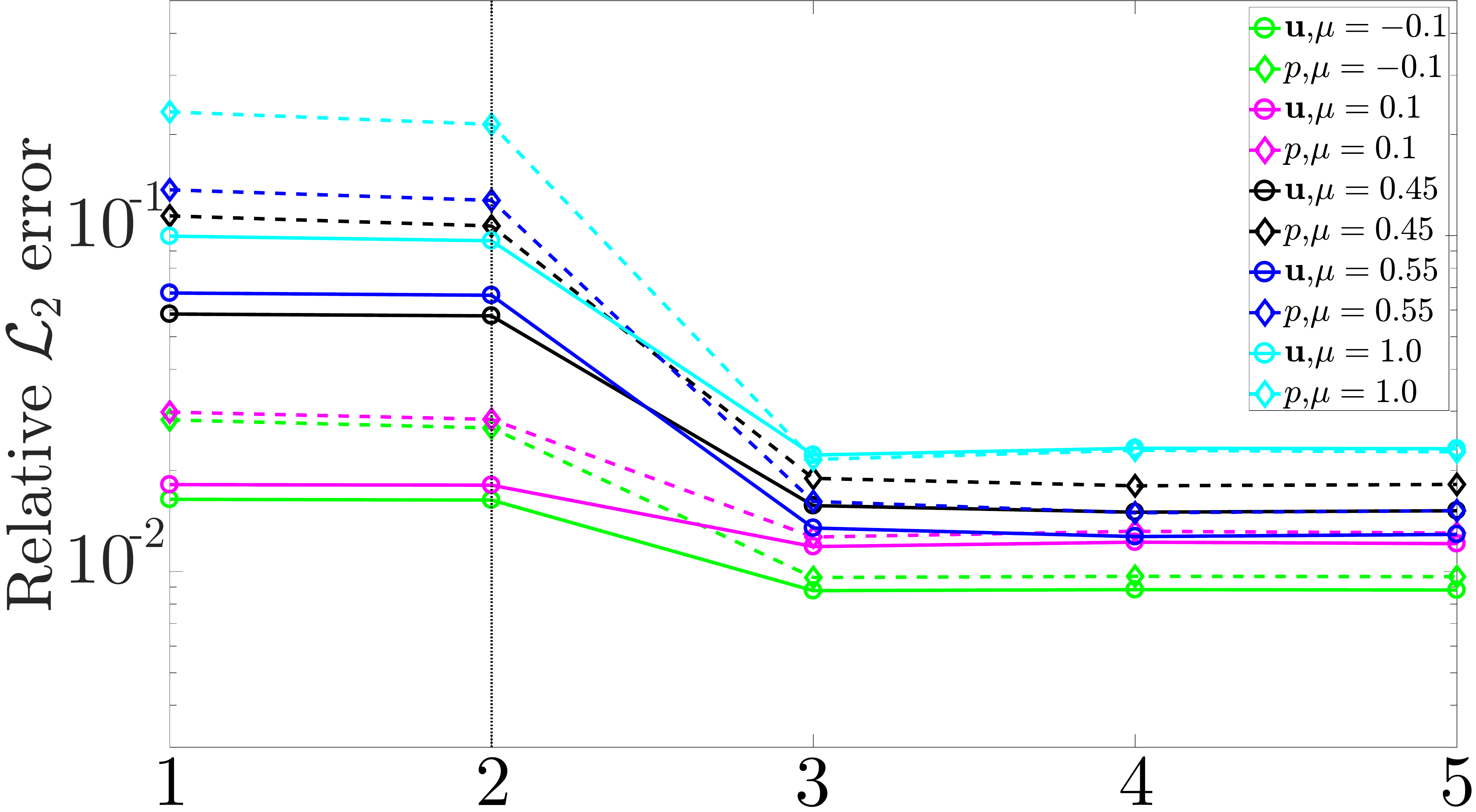}\label{fig:sBendErrorMode}}
	\caption{Internal flow in the S-Bend with parametrised jet velocity. (a) Relative amplitude of the computed modes $\fu^m$ (black), $\fp^m$ (blue) and the combined amplitude of $(\fu^m,\fp^m)$ according to Equation~\eqref{eq:relAmplitude}. (b) Relative $\eltwo(\Omega)$ error of the PGD approximations of pressure and velocity with respect to the full-order solutions for different values of $\mu$.}
	\label{fig:sBendAmplitude-ErrorMode}
\end{figure}

Setting a tolerance of $10^{-3}$, \texttt{pgdFoam} computes three modes before fulfilling the stopping criterion for $\eta_{(u,p)}$, see Equation~\eqref{eq:relAmplitude}, as displayed in Figure~\ref{fig:sBendAmplitude}.

The PGD approximation obtained using three computed modes is compared with the full-order solutions given by \texttt{simpleFoam} for the values $\mu{=}-0.1$, $\mu{=}0.45$ and $\mu{=}1$ of the parameter under analysis.
In Figure~\ref{fig:sBendErrorMode}, the relative $\eltwo(\Omega)$ error for these configurations is reported. The numerical experiments confirm that an accuracy of $10^{-2}$ is achieved using one computed mode additionnally to the two terms accounting for the boundary conditions. It is worth noting that the first computed mode is two orders of magnitude more relevant than the following ones (Fig.~\ref{fig:sBendAmplitude}). Thus, after one computed mode, the additional terms only introduce limited corrections to the existing PGD approximation.

\begin{figure}%[ht]
	\centering	
	\subfigure[$\ppgd, \mu{=}-0.1$]{\includegraphics[width=0.3\textwidth]{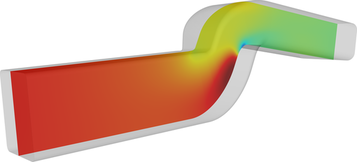}}
	\subfigure[$\ppgd, \mu{=}0.45$]{\includegraphics[width=0.3\textwidth]{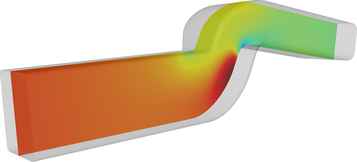}}
	\subfigure[$\ppgd, \mu{=}1$]{\includegraphics[width=0.3\textwidth]{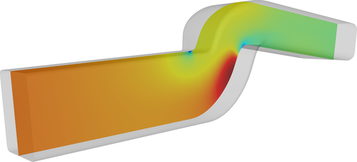}}

	\subfigure[$\pref, \mu{=}-0.1$]{\includegraphics[width=0.3\textwidth]{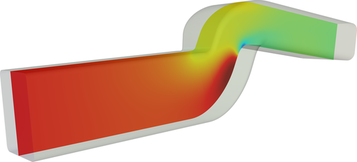}}
	\subfigure[$\pref, \mu{=}0.45$]{\includegraphics[width=0.3\textwidth]{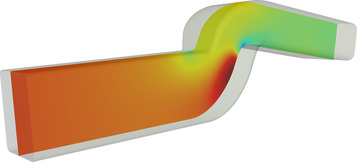}}
	\subfigure[$\pref, \mu{=}1$]{\includegraphics[width=0.3\textwidth]{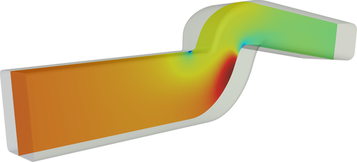}}
	
	\subfigure{\includegraphics[width=0.7\textwidth]{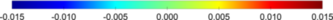}}
	\caption{Comparison of the PGD approximation (top) and the full-order solution (bottom) of the pressure field of the internal flow in the S-Bend with a jet configuration of $\mu=-0.1$, $\mu=0.45$ and $\mu=1$.}
\label{fig:sBend_Papproximation}
%\end{figure}
%
%\begin{figure}%[ht]
	\centering
	\subfigure[$\upgd, \mu{=}-0.1$]{\includegraphics[width=0.3\textwidth]{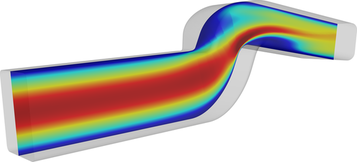}}
	\subfigure[$\upgd, \mu{=}0.45$]{\includegraphics[width=0.3\textwidth]{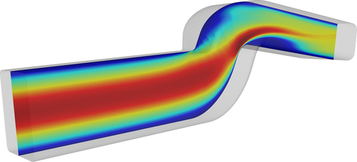}}
	\subfigure[$\upgd, \mu{=}1$]{\includegraphics[width=0.3\textwidth]{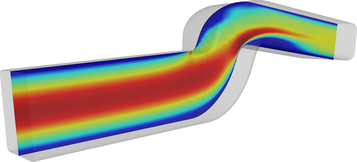}}

	\subfigure[$\uref, \mu{=}-0.1$]{\includegraphics[width=0.3\textwidth]{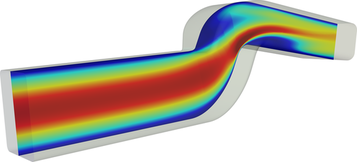}}
	\subfigure[$\uref, \mu{=}0.45$]{\includegraphics[width=0.3\textwidth]{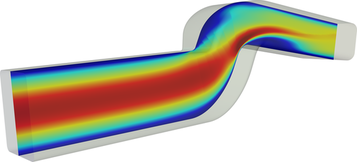}}
	\subfigure[$\uref, \mu{=}1$]{\includegraphics[width=0.3\textwidth]{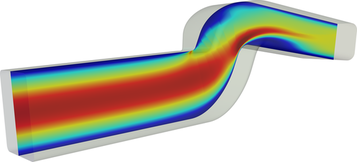}}
		
	\subfigure{\includegraphics[width=0.7\textwidth]{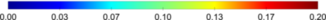}}
	\caption{Comparison of the PGD approximation (top) and the full-order solution (bottom) of the velocity field of the internal flow in the S-Bend with a jet configuration of $\mu=-0.1$, $\mu=0.45$ and $\mu=1$.}
\label{fig:sBend_Uapproximation}
\end{figure}

A qualitative comparison of the pressure and velocity fields computed using the PGD solution particularised for different values of the parameter $\mu$ and the corresponding full-order discretisations is presented in Figures~\ref{fig:sBend_Papproximation} and \ref{fig:sBend_Uapproximation}.

\begin{figure}%[ht]
\centering
	\subfigure[Error of the pressure drop]{\includegraphics[width=0.48\textwidth]{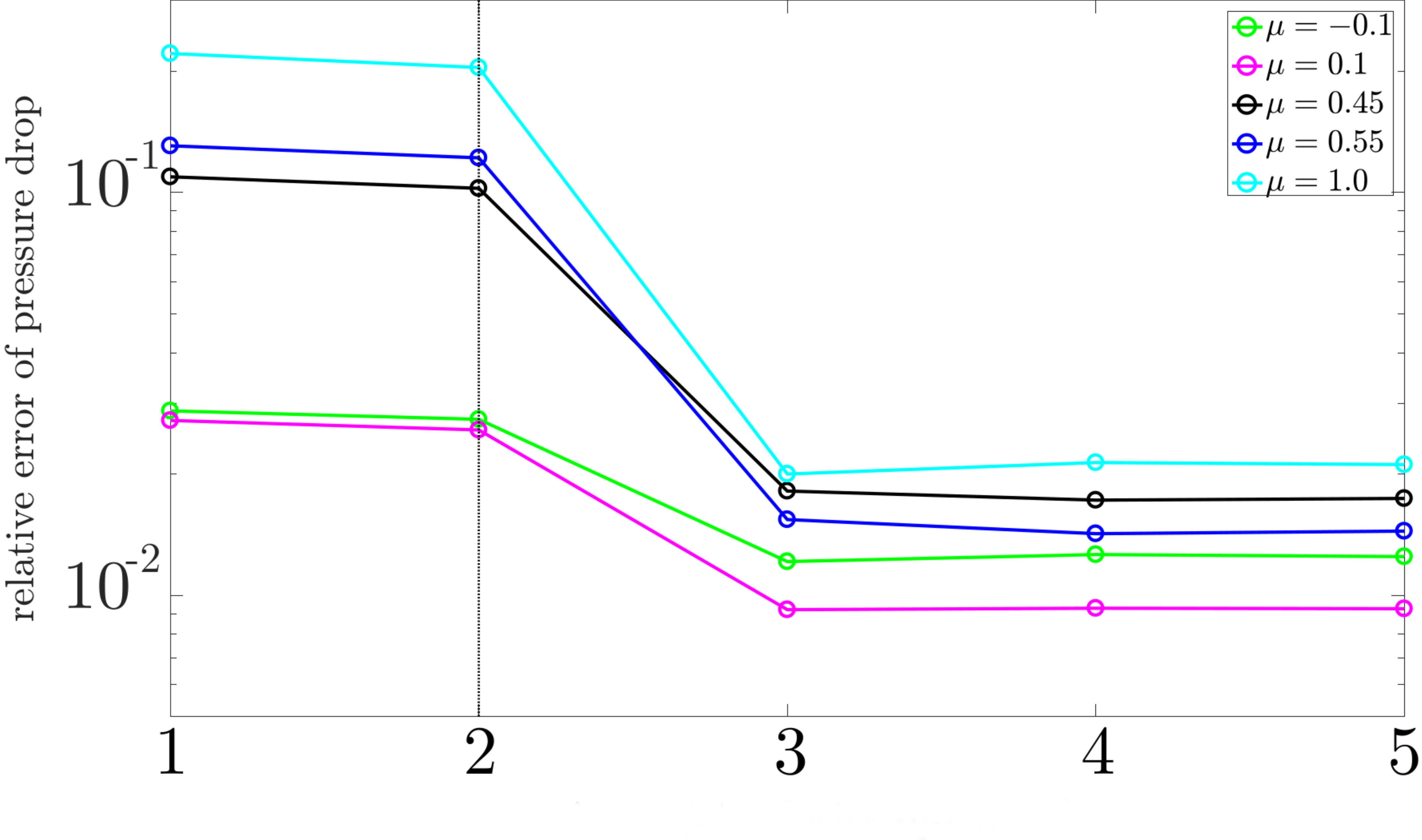}\label{fig:sBendError}}
	\subfigure[Pressure drop]{\includegraphics[width=0.48\textwidth]{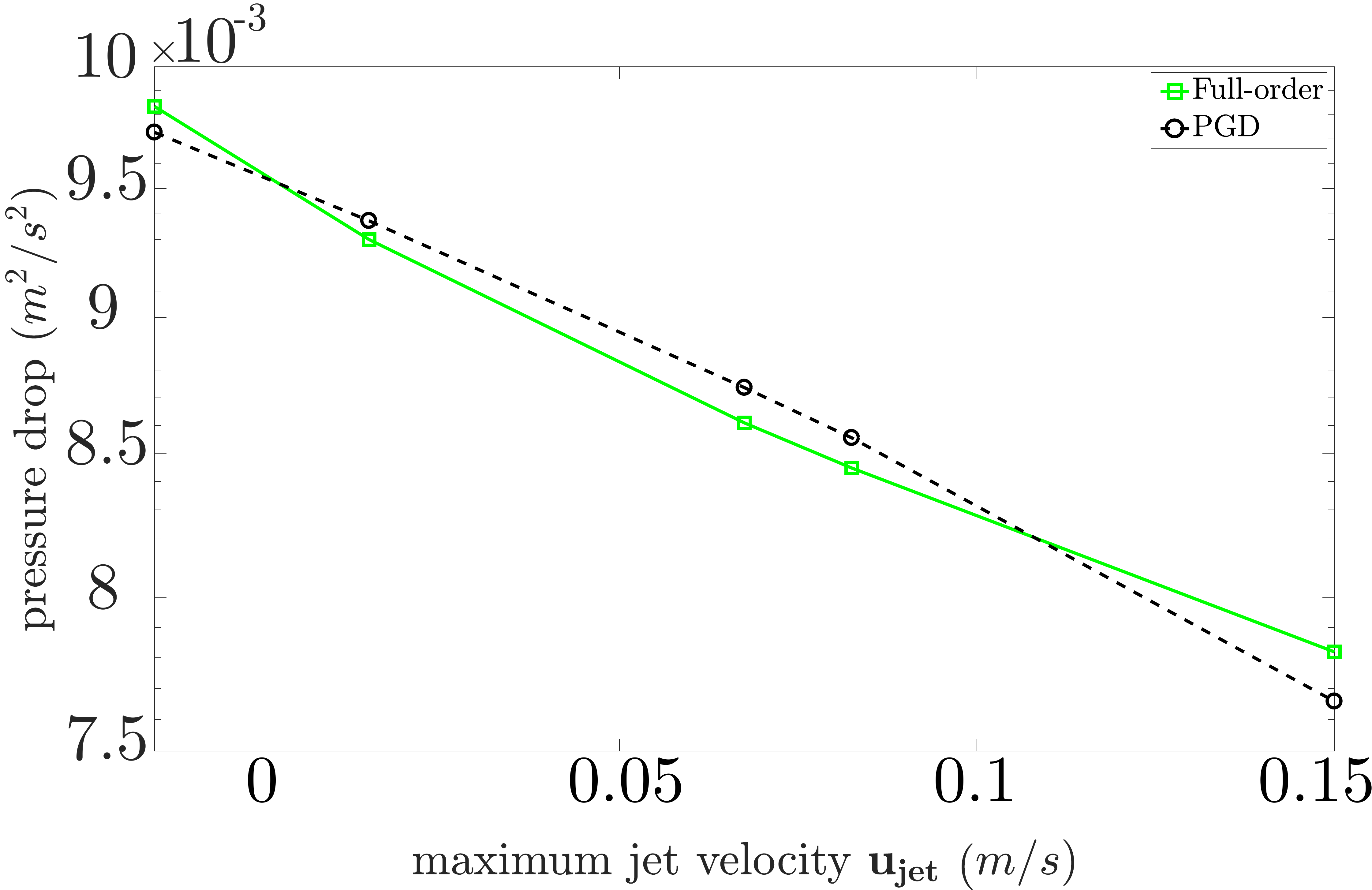}\label{fig:sBendDrop}}
	\caption{PGD approximation of the internal flow in the S-Bend with parametrised jet velocity. (a) Relative error of the pressure drop enriching the PGD modal approximation. (b) Pressure drop with respect to the maximum jet velocity.}
\label{fig:sBendError-Drop}
\end{figure}

As mentioned at the beginning of this section, the quantity of engineering interest in the analysis of this problem is the pressure drop computed along the duct.
The weighted average pressure drop is defined as
\begin{equation}\label{eq:sBendPdrop}
p_{\text{drop}}:=\frac{1}{A_{\text{in}}}\sum_{i=1}^{N_{\text{in}}}A_i p_i
\end{equation}
where $A_{\text{in}}$ is the area of the inlet surface, $N_{\text{in}}$ the number of faces $S_i, \, i{=}1,\ldots,N_{\text{in}}$ on the inlet patch and $p_i$, $A_i$ are the pressure and area on the face $S_i$, respectively.
For $\mu{=}-0.1$, $\mu{=}0.1$, $\mu{=}0.45$, $\mu{=}0.55$ and $\mu{=}1$, the pressure drop is evaluated as a particularisation of the generalised PGD solution and using the full-order solver \texttt{simpleFoam}.
Figure~\ref{fig:sBendError} presents the convergence history of the error in the pressure drop as a function of the number of modes in the PGD approximation.
It is straightforward to observe again that using the modes accounting for the boundary conditions and one computed mode is sufficient to capture the flow features of a wide range of parameters.
Moreover, by comparing the pressure drop with respect to the maximum velocity of the jet for different configurations with the corresponding values provided by the full-order solver, the capability of the discussed reduced-order strategy to accurately capture the evolution of a quantity of interest throughout the range of values of the parameter $\mu$ is confirmed (Fig.~\ref{fig:sBendDrop}).

%==========================================================================
\section{Concluding remarks}
\label{sc:Conclusion}
%==========================================================================

A nonintrusive PGD implementation in OpenFOAM has been proposed in the context of parametrised incompressible laminar flows.
The main novelty of such approach is represented by the seamless exploitation of OpenFOAM native SIMPLE solver, making the resulting reduced-order strategy suitable for application in a daily industrial environment.
The \texttt{pgdFoam} algorithm relies on the industrially-validated solver \texttt{simpleFoam} to compute the spatial modes of the solution, whereas the parametric ones are determined via the solution of a linear system of algebraic equations. 

The developed strategy has been validated using a manufactured solution to verify the optimal order of convergence of the PGD-FV approximation and a classical benchmark test case in the literature of CFD techniques for incompressible flows.
Moreover, the potential of the proposed PGD approach to rapidly and accurately simulate incompressible flows for different sets of user-defined parameters has been tested in the context of flow control problems. The \texttt{pgdFoam} algorithm has been applied both to an academic test case and an industrial one with a 3D geometrical model provided by Volkswagen AG. 

The proposed PGD methodology has proved to be able to compute an accurate reduced basis for the problems under analysis with no \emph{a priori} knowledge of the expected solutions. 
Moreover, it has shown robustness when dealing with a large range of values of the parameters, accuracy in capturing significant topological changes in the flow features and reliability in evaluating quantities of engineering interest, with an extremely reduced computing time.

%==========================================================================
\section*{Acknowledgements}
%==========================================================================
This work was partially supported by the European Union's Horizon 2020 research and innovation programme under the Marie Sk\l odowska-Curie Actions (Grant agreement No. 675919) that financed the Ph.D. fellowship of the first author.
The second, third and last author were also supported by the Spanish Ministry of Economy and Competitiveness (Grant agreement No. DPI2017-85139-C2-2-R). 
The second and last authors are grateful for the financial support provided by the Generalitat de Catalunya (Grant agreement No. 2017-SGR-1278).

%==========================================================================
%\bibliographystyle{plain}
\bibliographystyle{elsarticle-num}
\bibliography{Ref-PGD}
%==========================================================================

%==========================================================================
\appendix

%==========================================================================
\section{Separated representation of the residuals}
\label{sc:residuals}
%==========================================================================

Consider a separable expression of the source term $\bm{s}(\bx,\bmu) {:=} \eta(\bmu)\bm{S}(\bx)$.
For the spatial iteration, the residuals in separated form read as
\begin{subequations}\label{eq:spatialRes}
\begin{align}
&
%\begin{aligned}
%R_u^n
%:=& \alpha_4 \!\! \int_{V_i}{\bm{S} \, dV} 
%- \sum_{m=1}^n \sum_{q=1}^n \alpha_5^{mq} \!\! \int_{V_i}{\grad {\cdot} (\sigma_u^m\fu^m {\otimes} \sigma_u^q\fu^q) \, dV} \\
%&+ \sum_{m=1}^n \alpha_6^m \!\! \int_{V_i}{\grad {\cdot} (D \grad (\sigma_u^m\fu^m)) \, dV} 
%-  \sum_{m=1}^n \alpha_7^m \!\! \int_{V_i}{\grad (\sigma_p^m\fp^m)  \, dV} ,
%\end{aligned}
\begin{aligned}
R_u^n
:=& \int_{V_i}{\alpha_4 \bm{S} \, dV} 
- \sum_{m=1}^n \sum_{q=1}^n \alpha_5^{mq} \!\! \int_{V_i}{\grad {\cdot} (\sigma_u^m\fu^m {\otimes} \sigma_u^q\fu^q) \, dV} \\
&+ \!\! \int_{V_i}{\!\! \grad {\cdot} \Big( D \grad \Big(\sum_{m=1}^n \alpha_6^m\sigma_u^m\fu^m \Big) \! \Big) \, dV} 
-   \!\! \int_{V_i}{\!\! \grad \Big(\sum_{m=1}^n \alpha_7^m \sigma_p^m\fp^m \Big)  \, dV} ,
\end{aligned}
\label{eq:spatialResU} \\
&
R_p^n 
:= - \!\! \int_{V_i}{\!\! \grad {\cdot} \Big( \sum_{m=1}^n \alpha_7^m\sigma_u^m \fu^m \Big) \, dV} ,
\label{eq:spatialResP}
\end{align}
\end{subequations}
where the following expressions for the coefficients are devised
\begin{equation}\label{eq:coef-spatial-res}
\left\{\begin{aligned}
  \alpha_4 &:= \int_{\bI}{ \phi^n \eta  \, d\bI} ,&
  \alpha_5^{mq} &:= \int_{\bI}{ \phi^n \phi^m \phi^q \, d\bI} , \\
  \alpha_6^m &:= \int_{\bI}{ \phi^n \phi^m \psi \, d\bI} ,&
  \alpha_7^m &:= \int_{\bI}{ \phi^n \phi^m \, d\bI} .
  \end{aligned}\right.
\end{equation}

For the parametric iteration, the separated expression of the residuals is
\begin{subequations}\label{eq:paramRes}
\begin{align}
&
r_u^n := a_4 \eta + \sum_{m=1}^n \left( - \sum_{q=1}^n a_5^{mq} \phi^q + a_6^m \psi - a_7^m \right) \phi^m ,
\label{eq:paramResU} \\
&
r_p := - \sum_{m=1}^n a_8^m \phi^m ,
\label{eq:paramResP}
\end{align}
\end{subequations}
where the coefficients depend solely on the spatial modes, namely
\begin{equation}\label{eq:coef-param-res}
  \!\!\!\!\left\{\begin{aligned}
  a_4 &:= \int_{V_i}{\!\!\! \sigma_u^n\fu^n {\cdot} \bm{S} \, dV} , \\
  a_5^{mq} &:= \int_{V_i}{\!\!\! \sigma_u^n\fu^n {\cdot} \bigl[\grad {\cdot} (\sigma_u^m\fu^m {\otimes} \sigma_u^q\fu^q)\bigr] \, dV} , \\
  a_6^m &:= \int_{V_i}{\!\!\! \sigma_u^n\fu^n {\cdot} \bigl[\grad {\cdot} (D \grad (\sigma_u^m\fu^m))\bigr] \, dV} , \\
  a_7^m &:= \int_{V_i}{\!\!\! \sigma_u^n\fu^n {\cdot} \grad (\sigma_p^m\fp^m)  \, dV} , \\
 a_8^m &:=  \int_{V_i}{\!\!\! \sigma_p^n\fp^n \grad {\cdot} (\sigma_u^m\fu^m) \, dV} .
  \end{aligned}\right.
\end{equation}

%==========================================================================
\section{\texttt{simpleFoam}: the semi-implicit method for pressure linked equations in OpenFOAM}
\label{sc:SIMPLE}
%==========================================================================

In OpenFOAM, the steady Navier-Stokes equations are approximated by means of an iterative procedure, namely \texttt{simpleFoam}.
This algorithm implements the SIMPLE method proposed in~\cite{Patankar-PS:72}.
SIMPLE is a fractional-step Chorin-Temam projection method~\cite{Temam} that has been extensively studied in the literature~\cite{Guermond-Quartapelle:98a,Guermond-Quartapelle:98b}.
First, an intermediate velocity $\bu^k$ is computed starting from the momentum equation and neglecting the contribution of pressure, see Equation~\eqref{eq:SIMPLE-uInt}; second, the step involving the incompressibility constraint is rewritten in terms of a Poisson equation for the pressure $p$, see Equation~\eqref{eq:SIMPLE-PPE}; eventually, a correction is applied to the intermediate velocity field to determine the final value $\bu$ in Equation~\eqref{eq:SIMPLE-uFin}.
Special attention is required to impose the correct set of boundary conditions in each step of the algorithm~\cite{Laval-Quartapelle:90}.
\begin{subequations} \label{eq:SIMPLE}
\begin{align}
&
\left\{
\begin{aligned}
\frac{\bu^k-\bu^{k-1}}{\Delta t} 
+ \grad {\cdot} (\bu^k {\otimes} \bu^{k-1}) 
- \grad {\cdot} (\nu \grad \bu^k) 
&= \bm{s}       &&\text{in $\Omega$,}\\
\bu^k &= \bu_D  &&\text{on $\Gamma_D$,}\\
\bn {\cdot} (\nu \grad \bu^k) 
&= \bm{t}  &&\text{on $\Gamma_N$,}
\end{aligned}
\right.
\label{eq:SIMPLE-uInt} \\
&
\left\{
\begin{aligned}
\grad {\cdot} (\grad p) &= \frac{1}{\Delta t} \grad {\cdot} \bu^k       &&\text{in $\Omega$,}\\
\bn {\cdot} \grad p &= 0  &&\text{on $\Gamma_D$,} \\
\bn p &= \bm{0}  &&\text{on $\Gamma_N$,}
\end{aligned}
\right.
\label{eq:SIMPLE-PPE} \\
&
\bu = \bu^k - \Delta t \grad p .
\label{eq:SIMPLE-uFin}
\end{align}
\end{subequations}
Note that the algorithm in Equation~\eqref{eq:SIMPLE} may also be rewritten in the framework of an algebraic splitting method~\cite{Quarteroni-QSV:00}. For a complete introduction to the subject, interested readers are referred to~\cite{Donea-Huerta}.

\end{document}